\documentclass[11pt]{article}
\usepackage[english]{babel}

\usepackage{amsmath} 
\usepackage{tabularx}
\usepackage{amssymb}
\usepackage{natbib} 
\usepackage{amsfonts}
\usepackage{relsize}
\usepackage{color,soul}
\usepackage{array}
\usepackage{pict2e}
\usepackage{amsthm}
\usepackage[table]{xcolor}
\usepackage{multicol}
\usepackage{csvsimple}
\usepackage{array}
\usepackage{lscape}
\usepackage{tikz}
\usepackage{diagbox}
\usepackage{enumerate}
\usepackage{array}
\newcolumntype{P}[1]{>{\centering\arraybackslash}p{#1}}
\usepackage{graphicx}
\usepackage{array}
\usepackage{makecell}
\newcolumntype{x}[1]{>{\centering\arraybackslash}p{#1}}
\usepackage{tikz}

\usepackage{gensymb}

\usepackage{hyperref}
\usepackage{url}
\usepackage{enumerate}
\usepackage{graphicx}
\usepackage{epstopdf}
\usepackage{epstopdf} 

\usepackage{geometry}
\usepackage{graphicx,caption,scalerel}
\usepackage[nottoc]{tocbibind}
\usepackage{tikz}
\usetikzlibrary{shapes, shadows, arrows}
\tikzstyle{block}=[draw,rectangle,fill=blue!5,text width=12 em,text centered, minimum height=12mm, node distance=5 em]
\tikzstyle{line} = [draw,-latex']
\usetikzlibrary{calc,shapes.misc,shapes.geometric}

\usepackage{soul} 
\usepackage{xcolor} 
\sethlcolor{yellow} 
\usepackage{amsmath}
\usepackage{xcolor}

\begin{document}
	\title{\textbf{Exploration of  Hepatitis B Virus Infection Dynamics through an Intracellular  Model}}
	\author{Rupchand Sutradhar* and   Durga  Charan Dalal}
	\date{Department of Mathematics, Indian
		Institute of Technology Guwahati,
		Guwahati, Assam, India}
	\maketitle
	\noindent
	Corresponding author:
	Rupchand Sutradhar, Department of
	Mathematics, Indian Institute of
	Technology Guwahati, Guwahati, Assam,
	781039, India.
	Email: rsutradhar@iitg.ac.in\vspace{0.2cm}\\
	
	\section{Abstract}
	 Analysis of the cell population generally provides average information about  viral infection in a host whereas the intracellular model  captures the individual cellular responses. The primary goal of this study is to comprehensively analyze  the intracellular dynamics of hepatitis B virus (HBV) infection and  to identify the most influential factors.
	In this study, an intracellular HBV infection dynamics model is proposed by considering several intracellular steps that are observed in the virus life cycle. Upon comparison with the experimental data, it is seen that the model solutions  exhibit a good agreement.  The well-known fourth-order    \textit{Runge-Kutta}
	method is applied to numerically solve the proposed model. The effects of  HBx proteins, dslDNA-containing intermediates, intracellular delay and initial concentration of cccDNAs are explicitly studied. In order to identify the most positively and also the most negatively sensitive parameter of the proposed model, the global sensitivity analysis  is performed using  the widely-used method, Latin hypercube sampling-partial rank correlation coefficients.
As a result, it is observed that HBx proteins  have notable impacts on the dynamics of the infection, whereas intracellular delay and dslDNA-containing intermediates may not significantly affect the  infection. This study also suggests that sub-viral particles could potentially contribute to the progression of the infection. Furthermore, recycling of capsids (an intracellular process perhaps unique to the HBV life cycle, where a portion of the newly produced capsids return to the nucleus and amplify the cccDNAs) is found to play an important role in enhancing the infection.
	\vspace{0.5cm}\\
	Keywords: Intracellular model, Hepatitis B, cccDNA, Global sensitivity analysis, Partial rank correlation coefficient
	
%
%

	\section{Introduction}
	The non-cytopathic hepatitis B virus (HBV) belongs to the family of Hepadnaviridae.
	It causes various kind of serious liver diseases, such as cirrhosis,  hepatocellular carcinoma (HCC) or liver cancer. HBV infection is one of the most common liver infections and a global public health threats. It has the potential to spread 50 to 100 times more than HIV/AIDS \citep{WHO_2021}. Throughout the world, nearly two billion people have been affected by HBV. Every year  around 1.2 million people become newly infected despite the existence of effective prophylactic vaccine. Almost 254 million people are chronically infected and nearly 820,000 people die each year due to this viral infection and accompanying consequences, such as liver cancer \citep{WHO_2021}. At present, two distinct therapeutic strategies are used for the treatment of chronic hepatitis B (CHB) patients. These include  standard interferon-alpha (IFN-$\alpha$) or pegylated interferon-alpha (Peg-IFN-$\alpha$) and direct-acting antivirals composed of nucleoside analogs (NAs), such as lamivudine (LAM), telbivudine (LDT), entecavir (ETV), adefovir dipivoxil (ADV) and tenofovir disoproxil fumarate (TDF) \citep{2019_Zoulim}.  There are many limitations to the current treatment options, which often fail to provide long-term virologic control \citep{2008_Papatheodoridis}. The majority of the infected people require lifelong therapy and   remain stigmatized with HBV infection. Discontinuation of  the treatment will lead to a rebound of the virus \citep{2013_Chevaliez}. Stability and persistence of covalently closed circular DNAs (cccDNAs) that reside in the nuclei of infected cells is one of the major obstacles to prevent  this viral infection \citep{2005_zoulim_new}.

	Mathematical models  help understand the underlying mechanisms of disease progression and transmission. Most of the previous studies  on this viral infection \citep{1996_Nowak,2006_Murray,2018_Chenar,2021_Ashish,2023_Sutradhar_Fracational,2008_Wang,2024_Sutradhar,2014_Qiang}  mainly focused on the cell population. 
    Recently, \citet{2024_Ciupe_mice} developed a couple of  within-host models of HBV infection to predict the immune mechanisms responsible for the reduction in serum HBV and HBsAg levels.
    By analyzing the cell population, one can get a general overview about the viral infection whereas intracellular dynamics provides how individual cell responds to the viral infection \citep{2020_Brandt}. 
	In a cell population, there are various types of cells that  differ in identity, state, function, etc. as cells are heterogeneous. The variation in DNA sequence  is one of the reasons for this heterogeneity \citep{2010_raj_variability}. 
	In the literature, it is observed  that when  cell population is considered to study any kind of viral infections, the total number of cells is  generally divided into two classes: uninfected and infected cells. Due to this classification, intracellular dynamics of the infection can not be explained by those models \citep{1996_Nowak,2014_ciupe_antibody, 2018_Chenar,2017_Goyal,2024_Sutradhar_On} that are  designed based on the consideration of  cell populations because these models assume that all infected cells share the same characteristics. By modeling intracellular processes, one may study (i) the effects of each component of the virus and  parameters that are involved in the intracellular dynamics, (ii) the efficiencies and modes of action of antiviral therapies, and (iii) the influence of host factors (host DNA repair mechanism, immune system, signaling pathways, etc.) utilized by the virions throughout the  stages of the viral life cycle, etc \citep{2024_Ciupe}. Moreover, the consideration of  intracellular mechanisms also helps in shaping treatment strategies and guiding drug development.
	Besides, it is demonstrated that if the intracellular interactions are neglected, it may lead to an over-estimation of viral lifespan and results in inaccurate assumptions about drug efficacy and mechanisms of action \citep{2013_Guedj,2013_Rong}. In addition,  it is important to note that  the intracellular dynamics models are essential to evaluate the efficacies of those drugs that  target specific stages of the viral life cycle  \citep{2013_Guedj,2019_Cardozo}. 
%
	In 1940,  \citet{1945_delbruck_burst} considered the phage-infected E.coli cells to study the heterogeneity in virus infected cells. In that experiment, it was shown that the amount of progeny viruses produced from each infected cell  differs significantly.
	According to the studies on influenza A virus \citep{2015_heldt_single}, foot-and-mouth disease virus (FMDV) \citep{2018_xin_single}, and poliovirus \citep{2014_schulte_single}, virion levels show considerable differences between individual cells. 
	
	Nowadays, intracellular dynamics  analysis  become a significant milestone in many fields, particularly in virology.
	The followings are some key advantages of the  intracellular dynamics analysis in virology:
	\begin{enumerate}[(i)]
		\item  Bulk cell population analysis may overlook  responses of individual cell, whereas intracellular dynamics analysis provides a detailed view of reaction of each cell at the micro level.  
		\item The  intracellular components of the virus having  significant influence  can be identified.
		\item One can determine the most sensitive parameter  involved in the infection.
	\end{enumerate}

In case of HBV infection, some aspects of intracellular dynamics have been modeled previously \citep{2011_Jun,2016_jun_nakabayashi}. In a recent study, \citet{2015_Murray} proposed an intracellular  model to describe HBV replication in individual hepatocyte and  examined cccDNA formation through a stochastic agent-based approach using ordinary differential equations. \citet{2016_Goyal_cell_to_cell} has been extended the model to study the impacts of cell-to-cell transmission (CCT) route during infection. They  reported that while CCT has no significant impact on the establishment of infection, it  impedes its clearance and  also contributes to the development of fulminant hepatitis. The model presented by \citet{2015_Murray} was further extended to a multi-scale  hepatocyte sinusoid mathematical model to  investigate the impacts of complex intrinsic spatial heterogeneity of the liver on the HBV infection
dynamics \citep{2017_Cangelosi}.  \citet{2021_Ali} conducted global sensitivity analysis and identified the key parameters of the Murray and Goyal's model \citep{2015_Murray} by using two different techniques: (i) sampling-based sensitivity analysis method: Latin hypercube sampling-partial rank correlation coefficient (LHS-PRCC) \citep{2008_marino_methodology}, and (ii) variance-based Sobol method \citep{2001_Sobol}. In 2018, \citet{2018_Goyal_integration} demonstrated through an intracellular model that during occult HBV infection, the CCT is a primary mechanism for driving the spread of the infection.
Recently, using single-cell laser capture micro dissection (scLCM),  \citet{2020_Balagopal_single,2020_Balagopal_HBV_HIV}  informed  that the cccDNA transcription was significantly reduced or even silenced in those individuals undergoing long-term NUC therapy. Moreover, the computational analysis of inferred interactions suggested that an altered adaptive immune response associated to CD4 cytotoxic T lymphocytes emerges in patients who are functionally cured \citep{2024_Narmada_Single-cell landscape}.  Therefore, intracellular models have wide applications in understanding the viral dynamics.

	In this study,  based on the available biological and clinical  findings about the HBV intracellular life cycle  \citep{2021_tu_hepatitis,2020_Alexander,2017_Allweiss,2011_lucifora_hepatitis,2017_Thomas_hbv,2018_Chunkyu_hepatitis}, an improved intracellular dynamics model is proposed.  In the previous studies \citep{2011_Jun,2015_Murray,2016_jun_nakabayashi,PPR:PPR672105}, it is observed that the following intracellular  processes 
  
    \begin{itemize}
        \item The production of  rcDNA-containing capsids from viruses,
        \item Transcription of cccDNAs into viral mRNAs,
        \item Reverse-transcription of pgRNA-containing capsids results in single-stranded DNA capsids,
        \item The conversion of single-stranded DNA-containing capsids to double-stranded DNA-containing capsids,
        \item Translation of 2.4 kb and 2.1 kb mRNAs into surface proteins,
        \item The production of complete virus through the envelopment of double-stranded DNA-containing capsids by surface proteins
    \end{itemize}
   were considered as linear.  Due to the unavailability of experimental data and the lack of  sufficient supporting documents that suggest nonlinear behavior, in this study,  the  intracellular processes, including those mentioned above,  are modeled  under the assumption of linearity that offers a simpler and easier way to represent  the complex intracellular processes. 
In addition, the effects of many biological factors like  enzymatic reactions, interleukins, the packaging signal, microtubule networks,  importin-$\alpha$ and importin-$\beta$, etc. are also not considered. 
%
%
%
	 The model includes the roles of HBx proteins, the impacts of capsid recycling and the formation of double-stranded linear DNA (dslDNA)-containing capsids, etc. This model also considers the assembly of both non-infectious sub-viral particles (SVPs) and infectious virion particles. In addition, superinfection is  taken into account.  The proposed model is further extended by incorporating intracellular delay. 	The comparisons of model solutions with the experimental data of four mice show that the model solutions align closely with the experimental data.  
 Moreover, the effects of initial concentration of cccDNAs are   examined. One of the main objectives of this study is to identify the  most positively as well as most negatively sensitive parameter for each compartment of the model through global sensitivity analysis. This enables us to explore synergistic effects of  parameters and to suggest novel forms of monotherapy and combination therapy.

	
	%
	\section{Model formulation: materials and method}\label{model_formulation}
			\subsection{Intracellular dynamics model}\label{Label:Intracellular dynamics model}
			HBV is a member of the hepadnaviridae family and by virtue of its exceptional characteristics, it replicates through RNA intermediates in a similar way as it happens in case of retroviruses. In this way, the replication cycle of HBV is able to perpetuate the infection in hepatocytes through its unique features.
			HBV replication begins when the virus enters  the hepatocytes through the sodium taurocholate cotransporting polypeptide receptor by receptor-mediated endocytosis \citep{2012_Yan}. 
			Although, initially, HBV binds to heparan sulfate proteoglycans (HSPGs) with a low affinity \citep{2007_Schulze}. 
			This model does not explicitly include the roles of NTCP and HSPGs receptors in HBV entry. 
			Inside the hepatocytes, virus releases its core particles, \textit{i.e.}, relaxed circular DNA (rcDNA)-containing capsids.
			The number of viruses produced from a single infected hepatocyte and the quantity of rcDNA-containing capsids generated from incoming viruses from extracellular space  are designated by $V$ and $R$, respectively. It is considered that the viruses uncoat their core particles with the release rate $\alpha_1$. The natural  decay rate of rcDNA-containing capsids is denoted by $\delta_r$. The parameter $\alpha_2$ represents the rate at which rcDNAs are converted into   cccDNAs \citep{2007_Haitao}.
			Mathematically, this biological phenomenon can be expressed by the following differential equation:  
			\begin{align} \label{eq_R}
				\frac{dR}{dt}	= \alpha_1 V-\alpha_2 R-\delta_r R.
			\end{align}
			In the next step, in order to release the viral genome, HBV nucleocapsids  travel to the nucleus of hepatocytes. To overcome the high viscosity of the cytoplasm, HBV utilizes the microtubular network for efficient nuclear delivery \citep{2005_Dohner}.
			After entering the genome into the nucleus,  the rcDNAs are repaired through the host DNA repair mechanism \citep{2020_Alexander}  and  are converted into cccDNAs which are denoted by $C$. There are several  evidences which show that the persistence of cccDNAs is one of the major obstacles to prevent  this viral infection \citep{2017_Allweiss,2005_zoulim_new}. It is also believed that  cccDNAs have significant contribution in the  persistence of HBV infection as these are resistant to degradation and remain in the nuclei of the infected hepatocytes even after treatment is completed \citep{2019_Zhu}. 
			 Despite serological evidence of viral clearance, cccDNAs persist in individuals, contributing to the potential for viral reactivation \citep{2004_Werle-Lapostolle}.  These can also remain inside the hepatocytes for months or even years. There are different ways in which cccDNAs are lost, such as cell proliferation, cell death due to cytolytic immune response, cell cure due to non-cytolytic immune  response, and natural death of infected cells \citep{2021_Lythgoe}. In this study, only the natural decay  of cccDNAs $(\delta_c)$ is taken into account.
			This intracellular mechanism can be formulated as follows:
			\begin{align} \label{eq_C}
				\frac{dC}{dt}	= \alpha_2 R+k_1e^{-\lambda S_p} D-\delta_c C.
			\end{align}
			Here, the  term $k_1e^{-\lambda S_p} D$ represents the recycling of rcDNA-containing capsids. The details about the recycling of capsids are discussed in the formulation of equation \eqref{eq_D}. In the ensuing step,  cccDNAs  use the cellular transcriptional tools to produce viral RNAs including pregenomic RNA (pgRNA) and other messenger RNAs (mRNAs) which are the key components in the production of viral proteins \citep{2009_Levrero}.  In general, three  subgenomic mRNAs (0.7 kb mRNA, 2.1 kb mRNA and 2.4 kb mRNA) and two gemomics mRNAs of 3.5 kb are transcribed from cccDNAs \citep{2007_Block}. These mRNAs are heterogeneous and positively oriented \citep{2016_lamontagne_hepatitis}. The incorporation of these viral mRNAs in the model is crucial due to their pivotal roles in the viral life cycle, including the involvement of surface proteins (HBsAg) production, 
			the activation of cccDNA through HBV X proteins (HBx), etc \citep{2015_Nassal}. For ease-of-use in model formulation, we denote 3.5 kb mRNA by $R_g$. As both  2.4 kb mRNA and  2.1 kb mRNA produce surface proteins, these are treated as a single compartment by $R_s$. Similarly, 0.7 kb mRNA is represented by $R_h$. The parameters $\lambda_{rg}, \lambda_{rs} ~\mbox{and}~ \lambda_{rh}$ reflect the  transcription rates of 3.5 kb mRNA, 2.1 \& 2.4 kb mRNA, and 0.7 kb mRNA, respectively.  HBx  proteins that are denoted by $H$ are produced from 0.7 kb mRNA and  prevent the cccDNA from becoming silent.
			HBx can also enhance the transcription rate of cccDNA \citep{2016_decorsiere_hepatitis}. Moreover, HBx inhibits the development of immune response to HBV infection, thereby preventing apoptosis of infected hepatocytes \citep{2011_lucifora_hepatitis,2014_feitelson_roles}. Therefore, HBx plays some important roles in HBV replication.  Keeping everything in mind, these complex replication steps  are demonstrated by the following system of equations:
			\begin{align} 
				&\dfrac{dR_g}{dt}= \lambda_{rg} \Phi C-\mu_1R_g P-\delta_{r_g}R_g,\label{eq_R_g}\\
				&\dfrac{dR_s}{dt}= \lambda_{rs} \Phi C+\lambda_{sdl} D_L-\lambda_{s_p} R_s-\delta_{r_s}R_s,\label{eq_R_s}\\
				&\dfrac{dR_h}{dt}= \lambda_{rh} \Phi C-\delta_{r_h}R_h,
			\end{align}
			where $\delta_{r_g},~\delta_{r_s}~\mbox{and}~\delta_{r_h}$ indicate corresponding natural decay rates of $R_g,~ R_s~\mbox{and}~R_h$. Here,
			$\Phi$ denotes the volume fraction of active cccDNAs. The de-silencing of cccDNAs depends on the concentration of HBx proteins \citep{2015_Nassal}. Recently,  \citet{2020_Fatehi} modeled the silencing and de-silencing mechanisms of cccDNAs  in the presence of HBx proteins by the following equations:
					$$\text{(A1) Silencing:~~~~~~} \text{cccDNA} \xrightarrow{\dfrac{\lambda_{off}}{1+s\bar{X}}} \text{silenced~cccDNA},$$
				$$\text{(A2) De-silencing:~~~} \text{Silenced~cccDNA} \xrightarrow{\lambda_{on}(1+s\bar{X})} \text{cccDNA},$$			
			where $\lambda_{off}$ and $\lambda_{on}$ denote the silencing and de-silencing rates of cccDNAs, respectively, and $s$ indicates
			the efficiency of the available numbers of HBx proteins according to the model proposed by  \citet{2020_Fatehi}. Following equations (A1) and (A2),  it is considered that 
			$$\Phi=\left(1-\dfrac{1}{\dfrac{1}{1-\Phi_{\mathsmaller{0}}}+H}\right)=\dfrac{\Phi_0+H(1-\Phi_0)}{1+H(1-\Phi_0)}~ \text{with}~  \Phi_0~ \text{being the initial volume fraction of active cccDNAs}.$$ 
			When the concentration of HBx proteins increases, the value of $\Phi$ also increases. If $H=0$, $\Phi$ becomes $\Phi_0$. The detailed derivation of $\Phi$ is provided in Section 2 of the supplementary material. The dslDNAs which are denoted by $D_L$ can also produce surface proteins but may not be able to produce functional pgRNAs due to some mutations that are introduced when it is converted into cccDNAs \citep{2017_Thomas_hbv}. cccDNAs which are produced from dslDNA are functionally defective, as non-homologous end-joining pathway is
			error-prone \citep{2015_Nassal}.  It is assumed that  2.1 and 2.4 kb mRNAs are produced from  dslDNAs with production rate $\lambda_{sdl}$.
			The meaning of two terms $\mu_1R_g P$ in equation \eqref{eq_R_g}
			and $\lambda_{s_p} R_s$ in equation \eqref{eq_R_s} are explained later.
			
			As a result of translation of these mRNAs by ribosomes, viral proteins  are synthesized \citep{2020_Fatehi}. A portion of 3.5 kb mRNAs is translated into viral polymerases ($P$), while another portion is translated into core proteins ($C_p$) \citep{2011_Lewellyn}. $\lambda_p$  and $\lambda_c$ indicate the subsequent translation rates of the polymerases and core proteins. It is considered that  0.7 kb mRNAs is translated into HBx proteins  with translation rate $\lambda_{h}$ and  degrade naturally at rate $\delta_h$.
			The corresponding governing equations  for polymerases, core proteins and HBx proteins  are formulated as 
			\begin{align} \label{eq_P,Cp}
				~~&\frac{dP}{dt}	= \lambda_p R_g-\mu_1R_gP-\delta_p P,\\
				~~&\frac{dC_p}{dt}=\lambda_c R_g-\mu_2m ZC_p-\delta_{c_p}C_p,\\
				~~&\frac{dH}{dt}= \lambda_h R_h-\delta_h H,
			\end{align}
			where, $\delta_p$ and $\delta_{c_p}$ are the decay rates of polymerase and core protein, respectively.  After translation,  3.5 kb mRNAs are reverse-transcribed to viral genome DNA by viral polymerases. The pgRNA and polymerase form a 1:1 ribonucleoproteins which  are generally termed as RNP complexes \citep{2017_Hu}. These are assembly competent. These RNP complexes are denoted by $Z$ and  $\mu_1$ denotes the interaction  rate between  3.5 kb RNA and the polymerase. The mass balance equation for RNP complexes is given by
			\begin{align}
				&\frac{dZ}{dt}	=\mu_1R_gP-\mu_2 ZC_p-\delta_z Z.
			\end{align}
			Here, $\delta_z$ reflects the decay rate of RNP complex.  In the next step, RNP complexes are encapsidated by core proteins (HBcAg) to form  
			pgRNA-containing capsids (pgNCs) with interaction rate $\mu_2$. In this model, it is assumed that nucleocapsids contain $m=120$ core proteins following a $T = 4$ architecture.  We ignore the $T = 3$ particles in this study because  these particles occur in only 5\% of cases \citep{2017_patel, 2024_Yang_Structural}.
			 The pgNCs are also known as immature nucleocapsids and is denoted by the symbol $P_g$. 
			Consider that $\delta_{p_g}$ is the decay rate of  pgRNA-containing capsids. Based on the existing biological studies \citep{2017_Hu}, it is observed that a portion of pgRNA-containing capsids is enveloped by the surface proteins and  are secreted from the infected hepatocytes as non-infectious viral particles. It is not considered in this study.  Mathematically, pgRNA-containing capsids adhere to the following dynamical equation:
			\begin{align} \label{eq_Pg}
				&\frac{dP_g}{dt}=\mu_2 ZC_p-\beta_1 P_g-\delta_{p_g}P_g.
			\end{align}
			The reverse-transcription is one of the key steps in this virus life cycle. Through this process, the viral RNAs are converted into the viral DNA. The pgRNA acts as a template for DNA synthesis \citep{2015_Hu}. This step is associated with a series of events involving both the host and  virus factors. After encapsidation by the core proteins,  it is assumed that the viral polymerase reverse-transcribes the pgRNA with a rate $\beta_1$ into  single-stranded DNA (ssDNA).
			The ssDNA-containing capsids are designated by $S$ with a degradation rate $\delta_s$. It is assumed that the double-stranded DNA-containing capsids (dsDNA) and dslDNA-containing capsids are   produced from ssDNA with the  same rate $\beta_2$ .  
			According to the findings of  \citet{2017_Thomas_hbv}, $\sim$90\% of nucleocapsids generates rcDNA-containing capsids after reverse-transcription, while the remaining $\sim$10\% produce dslDNA-containing capsids.  In order to determine the relative contributions of different types of rcDNA-containing capsids on cccDNA production, this model distinguishes between  rcDNA ($R$) generated from upcoming virus from extracellular space and newly produced rcDNA \textit{i.e.} dsDNA-containing capsids, represented by $D$.   Newly produced rcDNA-containing capsids  can either gain an envelope of HBsAg by passing through the endoplasmic reticulum, pre-Golgi compartment, and  are released as infectious virion particles into the blood, or  can recycle back to the nucleus \citep{2018_Chunkyu_hepatitis}. In case of low level of surface proteins, rcDNA-containing capsids deliver its content to the nucleus and increase the pool of  cccDNAs \citep{2016_Jianming}.  This process is known as the `recycling of  capsids'.
			$\lambda^{-1}$ denotes the average  level of  surface proteins ($S_p$) \citep{2015_Murray}. Here, $k_1$ stands for the recycling rate of capsids. 
			The reaction equations for ssDNA, dsDNA and dslDNA-containing capsids can be given by the following system of equations \eqref{eq_S}-\eqref{eq_D_L}. 
			\begin{align}
				&\frac{dS}{dt}	=\beta_1 P_g-\beta_2 S-\delta_s S,\label{eq_S}\\
				&\frac{dD}{dt}	=0.9\beta_2 S-k_1e^{-\lambda S_p}D-k_2\left(1-e^{-\lambda S_p}\right)D S_p-\delta_d D,\label{eq_D}\\
				&\frac{dD_L}{dt}	=0.1 \beta_2 S-\lambda_{sdl} D_L-\delta_{d_L} D_L.\label{eq_D_L}
			\end{align}
			In equation \eqref{eq_C} as well as in equation \eqref{eq_D}, the term $ k_1e^{-\lambda S_p} D$ represents the recycling of capsids \citep{2015_Murray}. The parameters $\delta_d$ and $\delta_{d_L}$ are the natural decay rates of dsDNA and dslDNA-containing capsids, respectively.

			HBV surface protein (HBsAg) is one of the key components of this virus. HBsAg serves as a critical biomarker in the therapeutic management of HBV infection, with its sustained reduction or loss  associated with improved clinical outcomes and often is considered as a  indicator of functional cure. In this context, \citet{2021_Kadelka} constructed a mathematical model, integrating HBsAg with patient data to assess the time-dependent effectiveness of ARC-520 in inhibiting the production of both HBsAg and HBeAg. These proteins are produced by the translation of subgenomic RNAs (2.4 kb and 2.1 kb mRNA) by ribosomes. 2.4 kb mRNAs are translated into large surface proteins whereas translation of 2.1 kb of mRNAs leads to middle and small surface proteins \citep{2012_Datta}. For simplicity, these three surface proteins are referred to as one compartment and designated by $S_p$ in this model. The natural decay rate of the surface proteins is denoted by $\delta_{s_p}$. The parameter $\lambda_{s_p}$ is considered to be the mean value of production rate of surface proteins.
			During the infection, both complete viral particles and  SVPs are produced simultaneously and released into the bloodstream. SVPs are secreted from the hepatocytes in much larger amounts compared to infectious virions \citep{2020_Tsukuda}.  SVPs are predominantly composed of viral surface antigens (HBsAg).
			 Mainly two types of SVPs are produced: (i) spherical particles measuring 25 nm in diameter, and (ii) filamentous with a diameter of 22 nm, which
			 can differ in length \citep{2013_Gerlich}. Spherical SVPs are made of 48 S (small) proteins whereas filamentous SVPs have 1:1:4 protein stoichiometry  for L (large), M (medium), and S (small) surface proteins \citep{2009_Short}.
			    All of these SVPs are non-infectious in nature. We denote the average production rate of SVPs by $\eta_{sp}$.  Accordingly, the pertinent dynamical equation is outlined below as,
			\begin{align}
				&\frac{dS_p}{dt}=\lambda_{s_p} R_s-k_2\left(1-e^{-\lambda S_p}\right)D S_p-\eta_{sp}S_p-\delta_{s_p}S_p.
			\end{align}
			The matured capsids are released from the infected hepatocytes into the extracellular space or blood  as infectious Dane particles or complete virions. The virions exit  via the cell's secretory pathway by exocytosis \citep{2007_Lambert} and follows the given mathematical relationship:
			\begin{align} \label{eq_V}
				&\frac{dV}{dt}	= k_2\left(1-e^{-\lambda S_p}\right)D S_p-\delta_v V,
			\end{align}
			where $\delta_v$  and  $k_2$ characterize the death and release rates of viruses, respectively. The governing equation for SVPs can be given by
			\begin{align} \label{eq_Sv}
				&\frac{dS_v}{dt}	=\eta_{sp}S_p -\delta_{sv} S_v.
			\end{align} 	
			Here, it is considered that SVPs (denoted by $S_v$) are  cleared at rate $\delta_{sv}$.  

			%
			%
			\subsection{Full dynamics model}
			Based on the law of mass action, the temporal change of each component of virus is formulated. The following  system of equations \eqref{full model} describes the intracellular dynamics of the HBV infection with the non-negative initial conditions:
%
			\begin{equation} \label{full model}
				\left.
				\begin{split}
					&\frac{dR}{dt}	= \alpha_1 V-\alpha_2 R-\delta_r R,\hspace{5.3cm}
					\frac{dC}{dt}	= \alpha_2 R+k_1e^{-\lambda S_p} D-\delta_c C,\\
					&\frac{dR_g}{dt}= \lambda_{rg} \Phi C-\mu_1R_g P-\delta_{r_g}R_g,\hspace{3.9cm}
					\frac{dR_s}{dt}= \lambda_{rs} \Phi C+\lambda_{sdl} D_L-\lambda_{s_p}R_s-\delta_{r_s}R_s,\\
					&\frac{dR_h}{dt}= \lambda_{rh} \Phi C-\delta_{r_h}R_h, \hspace{5.6cm}
					\frac{dH}{dt}= \lambda_h R_h-\delta_h H,\\
					&\frac{dP}{dt}	= \lambda_p R_g-\mu_1R_gP-\delta_p P,\hspace{4.9cm}
					\frac{dZ}{dt}	=\mu_1R_gP-\mu_2 ZC_p-\delta_z Z,\\
					&\frac{dC_p}{dt}=\lambda_c R_g-\mu_2m ZC_p-\delta_{c_p}C_p,\hspace{4cm}
					\frac{dP_g}{dt}=\mu_2 ZC_p-\beta_1 P_g-\delta_{p_g}P_g,\\
					&\frac{dS_p}{dt}=\lambda_{s_p} R_s-k_2(1-e^{-\lambda S_p})D S_p-\eta_{sp}S_p-\delta_{s_p}S_p,\hspace{1.2cm}
					\frac{dS}{dt}	=\beta_1 P_g- \beta_2 S-\delta_sS,\\
					&\frac{dD_L}{dt}	=0.1 \beta_2 S-\lambda_{sdl} D_L-\delta_{d_L} D_L,\hspace{4cm}
					\frac{dV}{dt}	= k_2\left(1-e^{-\lambda S_p}\right)D S_p-\delta_v V,\\
					&\frac{dD}{dt}	=0.9\beta_2 S -k_1e^{-\lambda S_p} D-k_2\left(1-e^{-\lambda S_p}\right)D S_p-\delta_d D, \hspace{0.2cm}
					 \frac{dS_v}{dt}	=\eta_{sp}S_p -\delta_{sv}S_v.\\
				\end{split}
				\hspace{0.01cm}
				\right\} 
			\end{equation}
		  The system \eqref{full model} is solved numerically using the fourth-order highly-accurate  \textit{Runge–Kutta} method.  The description of all model variables  are summarized in  Table \ref{Table_Variables}.	All the model parameters  are non-negative. The baseline values, units and sources of the  parameters are provided in  Table \ref{Table_parameters}.						
%
%
			\begin{table}[h!]
				\small
				\caption{State variables and their biological meanings}
				\begin{center}
					\begin{tabular}{||c|| l|| }
						\hline
						Variables  & ~~~~~~~~~Descriptions  \\ [.5ex]
						\hline
						$R$   	&	Number of rcDNA-containing capsids\\
						$C$		&	Number of cccDNAs\\	
						$R_g$ 	&	Number of 3.5 kb pgRNAs\\
						$R_s$	&	Number of 2.4 and 2.1 kb mRNAs\\
						$R_h$	&	Number of 0.7 kb mRNAs\\
						$H$  	&	Number of HBx proteins\\
						$P$		&	Number of polymerases \\
						$Z$		&	Number of RNP complexes\\
						$C_p$	&	Number of core proteins\\
						$P_g$	&	Number of pgRNA-containing capsids\\
						$S_p$   &   Number of surface proteins\\
						$S$     &	Number of single-stranded DNA-containing capsids\\
						$D$		&	Number of newly produced double-stranded DNA-containing capsids\\
						$D_L$	&	Number of  double-stranded linear DNA-containing capsids\\
						$V$ 	&	Number of Viruses\\
						$S_v$	&	Number of SVPs\\
						\hline		
					\end{tabular}
				\end{center}
				\label{Table_Variables}
			\end{table}
				\begin{table}
	\small
	\caption{Descriptions of model parameters, including their sources and units. 
		The  values of  parameters are taken from literature \citep{2015_Murray,2021_Fatehi} with some slight symbolic changes. The letter `P' indicates parameters.}
	\begin{center}
		\begin{tabular}{||c|m{6.5cm}|m{4.2cm}|l||}
			\hline
			P  & ~~~~~~~~~Descriptions & Value  & Source  \\ [.5ex]
			\hline
			$\alpha_1$  	& Production rate of rcDNA-containing capsid from viruses	& $0.03/24~ hour^{-1}$ & \citep{2015_Murray}\\
			\hline
			$\alpha_2$   & Production rate of cccDNA from rcDNA-containing capsids & $\log(2)/24~hour^{-1}$ & \citep{2015_Murray}\\
			\hline
			$\lambda^{-1}$   & Average $S_p$ level & $100000$ molecules/cell& \citep{2015_Murray}\\
			\hline
			$k_1$	 	& Recycling rate of rcDNA-containing capsids	& $\log(2)/24~hour^{-1}$&\citep{2015_Murray}\\
			\hline	
			
			$k_2$		& Binding rate of surface protein with rcDNA-containing capsids & $\log(2)/24~mole^{-1}hour^{-1}$ & \citep{2015_Murray}\\
			\hline
			$\lambda_{rg}$ & Production rate of 3.5 kb pgRNA from cccDNA & 25.71 $hour•^{-1}$	& \citep{2021_Fatehi}\\
			\hline
			$\lambda_{rs}$ 	& Production rate of 2.4 and 2.1 kb mRNA from cccDNA	& 37.5 $hour^{-1}$& \citep{2021_Fatehi}\\
			\hline
			$\lambda_{rh}$	& Production rate of 0.7 kb mRNA from cccDNA & 128.57 $hour^{-1}$	&\citep{2021_Fatehi}\\
			\hline
			$\lambda_{h}$	& Production rate of HBx proteins from 0.7 kb mRNA 	& 116.88 $hour^{-1}$& \citep{2021_Fatehi}\\
			\hline
			$\lambda_p$	& Production rate of polymerase from pgRNA	& 22.5 $hour^{-1}$ &\citep{2021_Fatehi}\\
			\hline
			$\lambda_{c}$ 	& Production rate of core protein from pgRNA	& 25.71 $hour^{-1}$ & \citep{2021_Fatehi} \\
			\hline
			$\lambda_{s_p}$ 	& Production rate of surface protein from 2.4 and 2.1 kb mRNA	&208.33 $hour^{-1}$&\citep{2021_Fatehi} \\
			\hline
			$\lambda_{sdl}$	& Production rate of 2.4 and 2.1 kb mRNA from dslDNA	& Assume (=$\lambda_{s_p}$) &~ -\\
			\hline
			$\mu_1$	    & Interaction rate 3.5 kb pgRNA and polymerase  & 5 $molecule^{-1}hour^{-1}$&\citep{2021_Fatehi}\\
			\hline
			$\mu_2$		& Interaction rate RNP complex and core protein	& Assume (=$\mu_1$) & ~ -\\
			\hline
			$\beta_1$ 	& Production rate of pgRNA-containing capsids	& 2.08 $hour^{-1}$ &\citep{2015_Murray}\\
			\hline
			$\beta_2$    & Production rate of ssDNA-containing capsids and dslDNA containing capsids &$\ln(2)/24 ~hour^{-1}$& \citep{2015_Murray}\\
			\hline
			
			
			$\delta_r$ 	& Decay rate of rcDNA-containing capsids	& $\ln(2)/24~ hour^{-1}$& \citep{2010_Xu_interferons}\\
			\hline
			$\delta_c$  & Decay rate of cccDNAs 	& $0.016/24~ hour^{-1}$ & \citep{2021_lythgoe_estimating}\\
			\hline
			$\delta_{r_g}$   & Decay rate of 3.5 kb pgRNAs &  $ \ln(2)/5~hour^{-1}$& \citep{2021_Fatehi}\\
			\hline
			$\delta_{r_g}$   & Decay rate of 2.4 and 2.1 kb mRNAs &  $ \ln(2)/3~hour^{-1}$& \citep{2021_Fatehi}\\
			\hline
			$\delta_{r_h}$	& Decay rate of 0.7 kb mRNAs	&  $\ln(2)/3~hour^{-1}$&\citep{2021_Fatehi}\\
			\hline
			$\delta_{h}$	& Decay rate of HBx proteins	&  $ \ln(2)~hour^{-1}$ & \citep{2021_Fatehi}\\
			\hline
			$\delta_p$		& Decay rate of polymerases & $ \ln(2)~hour^{-1}$	& \citep{2021_Fatehi}\\
			\hline
			$\delta_z$ 	    & Decay rate of RNP complexes &   $0.06~hour^{-1}$ &\citep{2016_jun_nakabayashi}\\
			\hline
			$\delta_{c_p}$	& Decay rate of core proteins	&  $\ln(2)~hour^{-1}$&\citep{2021_Fatehi}\\
			\hline
			$\delta_d$	    & Decay rate of dsDNA-containing capsids &$0.053/24~hour^{-1}$& \citep{2006_Murray}\\
			\hline
			$\delta_{p_g}$	& Decay rate of pgRNA-containing capsids	& Assumed ($=\delta_d$)& ~ -\\
			\hline
			$\delta_{s_p}$	& Decay rate of surface proteins	&  $ \ln(2)~hour^{-1}$& \citep{2021_Fatehi}\\	
			\hline	
			$\delta_s$		& Decay rate of ssDNA-containing capsids	&Assumed ($=\delta_d$)& ~ -\\
			\hline
			
			$\delta_{d_L}$	    & Decay rate of dslDNA-containing capsids & Assumed ($=\delta_d$)&~ -\\
			\hline
			$\delta_v$		& Natural death rate of viruses	& $0.158~hour^{-1}  $ & \citep{2006_Murray}\\
			\hline
			$\delta_{sv}$		& Natural death rate of SVPs	& ($=\delta_v$)   & \citep{2014_ciupe_antibody}\\
			\hline
			$\eta_{sp}$		& Exit rate of SVPs	& $2.9~hour^{-1}$ & \citep{2021_Fatehi}\\
			\hline
			$\Phi_0$		& Initial volume fractions of active cccDNAs	& 0.2 (Assumed) & ~ - \\
			\hline
			
		\end{tabular}
	\end{center}
	\label{Table_parameters}
\end{table}

			\section{Model validation}
\begin{figure}[h!]
		\centering		\includegraphics[height=11cm,width=16cm]{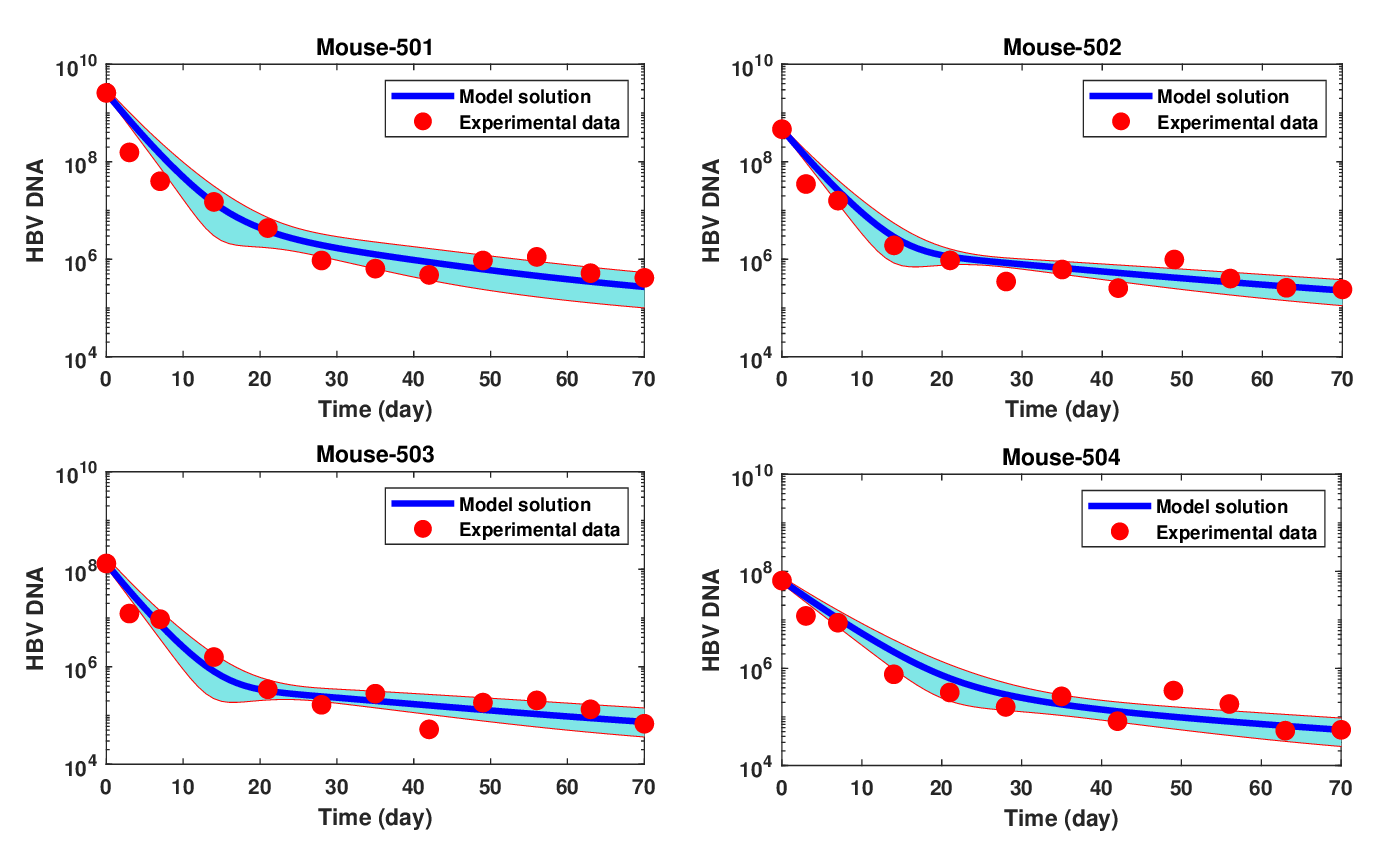}
		\caption{Calibration of the proposed model using experimental data collected from four humanized mice subjected to ETV treatment. 
			 Both model solutions and experimental data are presented on a logarithmic scale along the $y$-axis.  The shaded areas represent the 95\% confidence interval of the model solution.}
		\label{Experimental_Mouse}
	\end{figure}

			In order to validate the proposed model \eqref{full model}, it is extended by incorporating the effects of entecavir (ETV) that acts as a reverse-transcriptase inhibitor and can block the production of  dsDNA-containing capsids from pgRNA-containing capsids. The efficiency of ETV ($\epsilon$) satisfies $0<\epsilon\leq 1$. As a result of incorporation of ETV, the equations \eqref{eq_Pg} and  \eqref{eq_S} are modified as follows: 
			\begin{equation} 
				\left.
				\begin{split}		
					\dfrac{dP_g}{dt}&=\mu_2 ZC_p-(1-\epsilon)\beta_1 P_g-\delta_{p_g}P_g,\\
					\dfrac{dS}{dt}&=(1-\epsilon)\beta_1 P_g- \beta_2 S-\delta_s S.\\
				\end{split}
				\right\}
			\end{equation}
			Experimental data of four humanized mice are collected from the work of  \citet{PPR:PPR672105}. Each mouse was infected with HBV at $1.0\times 10^6$ copies.
			On the day 53 from post-inoculation, those mice displaying a sustained level of HBV in serum, were administered ETV continuously for 70 days. The treatment protocol involved daily dosing of 0.02 mg/kg ETV. The efficiency of ETV, as stated in the study by  \citet{PPR:PPR672105}, is recorded as 0.97, and this value is utilized in this study. By estimating the model parameters, we compare the model solutions with the experimental data from four humanized mice (labeled as Mouse-501, Mouse-502, Mouse-503, Mouse-504). The solutions as well as experimental data are plotted in Figure  \ref{Experimental_Mouse}. In order to enhance the clarity of the model validation, some statistical testings, such as $R-squared$, relative error, confidence intervals (CIs)  are performed. $R-squared$ is given by the following formula:
			$$R-squared = 1 - \dfrac{\displaystyle\sum_{i=1}^{N} \left(y_i^{\text{obs}} - {y}^{\text{pred}}\right)^2}{\displaystyle\sum_{i=1}^{N} \left(y_i^{\text{obs}} - \overline y_i^{\text{obs}}\right)^2},~\text{where}~\overline y_i^{\text{obs}}~\text{denotes the sample mean}.$$
		The values of $R-squared$ are shown in the second column of  Table \ref{Table: The values of R-squared}. All values  $R-squared$ are greater than 0.9. So, these $R-squared$'s values suggest a strong correlation between model solutions and experimental data. The relative errors are presented in third column of Table \ref{Table: The values of R-squared} for all mice. Additionally, we calculate 95\% CIs shown in Figure \ref{Experimental_Mouse}. The CIs are found to be  narrow, indicating precise parameter estimation.
			\begin{table}[h!]
				\caption{The values of $R-squared$.}
				\begin{center}
					\begin{tabular}{|c|c|c|}
						\hline
						Mouse     & The value of $R-squared$ & Relative error\\ \hline
						Mouse-501 & 0.9479 &    0.02093          \\ \hline
						Mouse-502 & 0.9557 &    0.01926          \\ \hline
						Mouse-503 & 0.9236 &   0.01740          \\ \hline
						Mouse-504 & 0.9185 &    0.02639         \\ \hline
					\end{tabular}
				\end{center}				
				\label{Table: The values of R-squared}
			\end{table}
		Based on the results of the statistical testing, it can be concluded that the model solutions agree well  with the  experimental data, and accurately capture the dynamics of the infection.

	\section{Results and Discussion}
			We have developed an intracellular model of HBV infection dynamics that incorporates essential steps of the viral life cycle.
			The onset of  infection is assumed to be  the initial condition which is represented by $$(R(0),C(0),R_g(0),R_s(0),R_h(0),H(0),P(0),Z(0),C_p(0),P_g(0),S_p(0),S(0),D_L(0),D(0),S_v(0),V(0))=$$
			$$(0,0,0,0,0,0,0,0,0,0,0,0,0,0,0,V(0)\neq0).$$ In order to study the effects of initial concentration of cccDNAs, we vary the initial concentration of cccDNAs. Otherwise, for the entire length of simulation, the stage of infection remains same. Even, same stage of infection is considered during global sensitivity analysis.
			According to \citet{2000_Nowak}, the lifespan of infected hepatocytes in  HBV infection that ranges from 10 to 100 days. 
			 In this study, we consider the maximum lifespan  of infected hepatocytes, \textit{i.e.}, 100 days to conduct the subsequent numerical experiments. 
			\subsection{Effects of initial concentration  of cccDNAs}\label{Effects of initial concentration  of cccDNA}
\begin{figure}[h!]
				\centering
				\includegraphics[height=14cm,width=18cm]{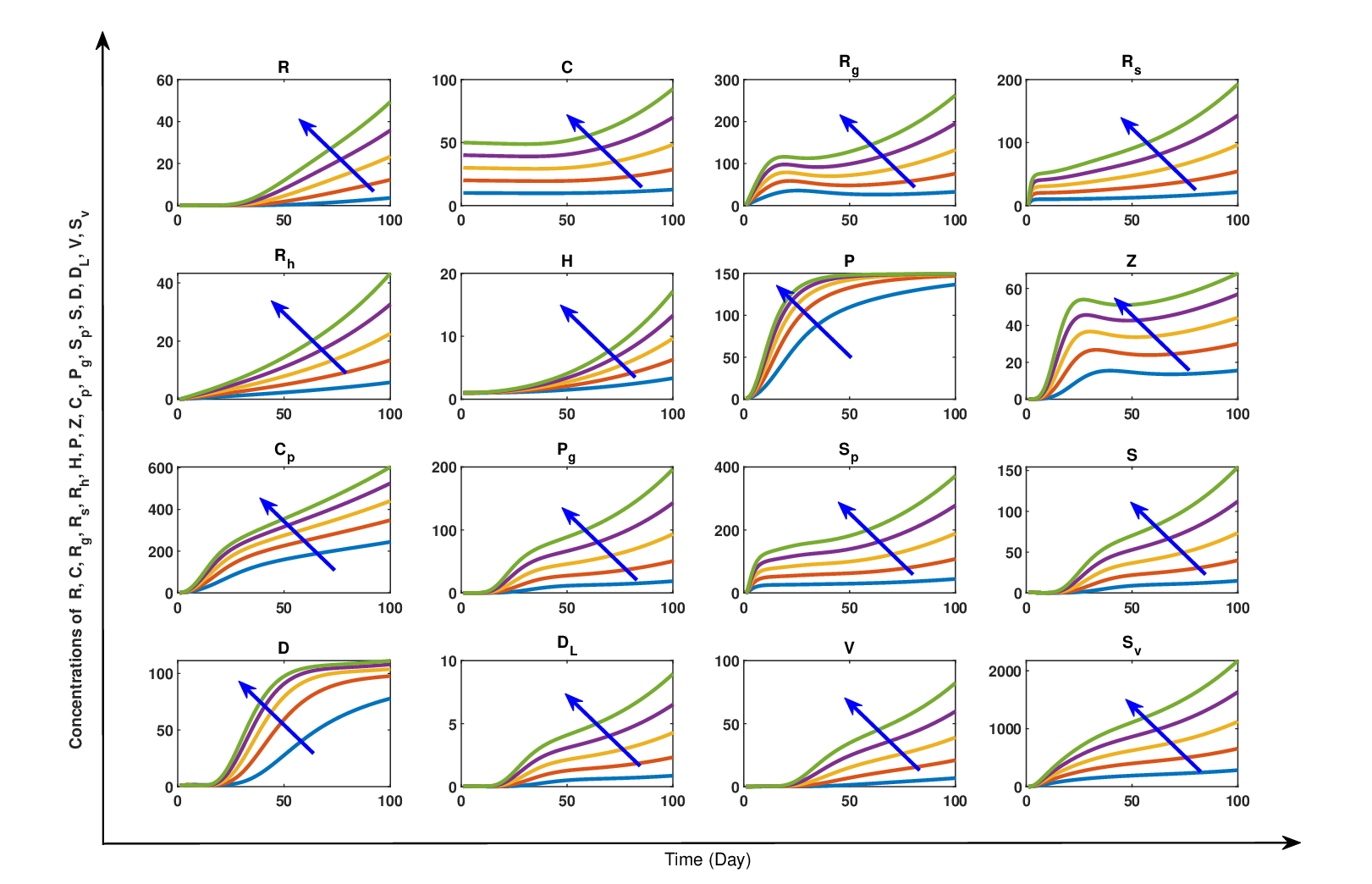}
				\caption{The effects of initial concentrations of cccDNAs on all components. Five different initial conditions $C(0)=10,20,30,40,50$ are considered here. Time and concentration of each component are taken  along $x-$axis and  along $y-$axis, respectively.  The blue arrows in each subplot indicate  direction of growth of the initial condition of cccDNAs.}
				\label{fig:The effects of initial concentration of cccDNA}
			\end{figure}

		The existence of cccDNAs is highly associated with the  persistence of infection in patients. It is thought to be a major factor for persistence of HBV infection as it is resistant to degradation and remains in the nucleus of infected cells even after treatment is completed. Due to its strong stability, cccDNAs are not lost in the time of cell division \citep{2004_locarnini_molecular,2007_glebe_recent}. Despite serological evidence of viral clearance, cccDNAs persist in individuals and act as a potential candidate for viral reactivation.  It can also remain inside the hepatocytes for months or even years \citep{2008_balsano_viral}. 
			The effects of initial concentration of cccDNAs on all compartments of the model \eqref{full model}  are demonstrated here.   The initial concentrations of cccDNAs  are varied keeping the initial concentrations of other compartments fixed. To this purpose, five distinct initial concentrations of cccDNAs  are considered. As a result, it is observed that the initial concentrations of cccDNAs significantly influence all compartments. The presence of a few copies of cccDNAs in the liver can re-initiate and blow-up the infection as shown in Figure \ref{fig:The effects of initial concentration of cccDNA}.
			\citet{2022_Chang} reported similar findings, as follows: even when just a single copy of cccDNAs remains in hepatocytes, HBV replication can lead to detectable viremia within a specific period of time following HBV reactivation. 
			 Though a small quantity of cccDNAs may not be clinically important, but  it remains  as a crucial factor in the persistence and transmission of HBV infection. 
			 Several  direct and indirect approaches, such as genome-editing approaches, strategies targeting host dependency factors or epigenetic gene regulation, nucleocapsid modulators and immune-mediated degradation   are  under  discovery, preclinical or early clinical development \citep{2020_Ligat}. Among gene-editing drugs, CRISPR/Cas9 has proven to be the most successful \citep{2016_Durantel}. 
%
			 Attacking HBV cccDNAs may be considered the holy grail for achieving a complete cure of hepatitis B \citep{2016_lucifora_attacking}.
			
			\subsection{Roles of HBx proteins  on HBV infection} \label{effect of HBx proteins}
\begin{figure}[h!]
			\centering
			\includegraphics[height=11cm,width=17cm]{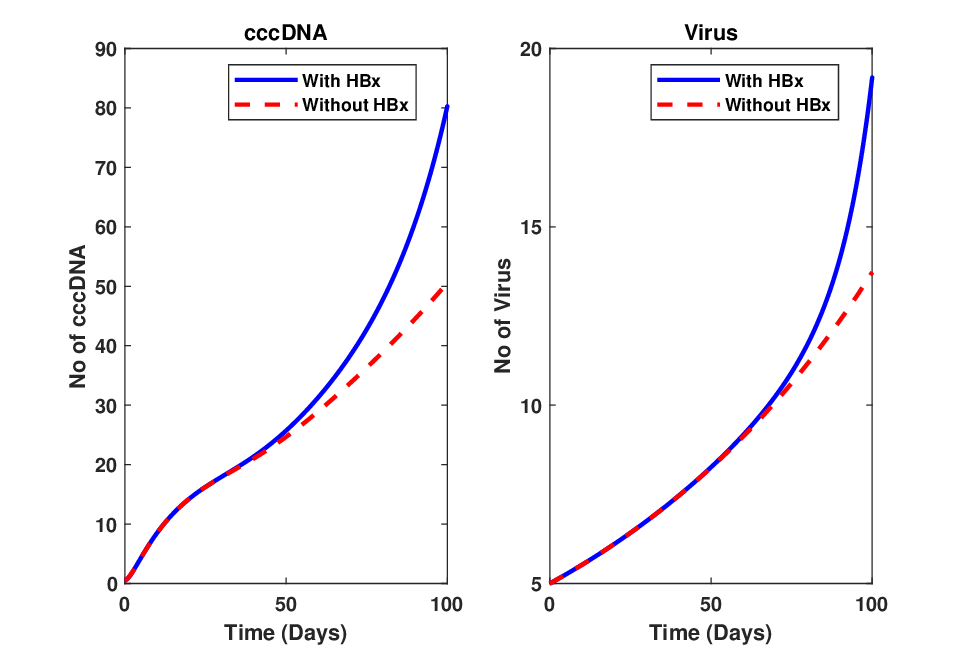}
			\caption{The effects of HBx proteins are demonstrated. Blue solid line  and red dotted line represent the solution of the system \eqref{full model} with and without considering the effects of the HBx proteins, respectively.}
		\label{The effects of HBx protein}
		\end{figure}

			HBx proteins contribute to the infection through various mechanisms, which have already been discussed in Section  \ref{Label:Intracellular dynamics model}. These proteins  can protect infected hepatocytes from immune-mediated apoptosis by suppressing  the immune systems in HBV infection \citep{2014_feitelson_roles}, although many aspects of HBx are still unknown.  Due to this, the impacts of HBx proteins are incorporated in the model \eqref{full model} to illustrate how it affects the intracellular components of the virus during infection.
			The simulation results show that the inclusion of HBx proteins substantially changes the concentration levels of  compartments of the model \eqref{full model}.
			  Figure \ref{The effects of HBx protein} illustrates the changes in the concentration of cccDNAs and viruses, while the variations in other components are presented in the supplementary material (Figure S1 and Figure S2).  In Figure \ref{The effects of HBx protein}, it is observed that the numbers of cccDNAs and viruses increase rapidly when HBx is taken into account. Similar  dynamics are also recorded for other compartments. In essence, the HBx proteins exert considerable influence on pathogenesis. Recently, \citet{2022_Sivasudhan_Hepatitis} also noted numerous pivotal functions of HBx in the persistent infection and advancement of HCC.
			Therefore,  considering the HBx protein as a therapeutic target may provide a promising avenue for future management of HBV infection. 
			 In this regard, the antiviral agent SC75741  is identified as a good inhibitor for  HBx proteins \citep{2024_Hamadalnil}.   Moreover, targeting the HBx proteins to maintain the silencing of cccDNA can be suggested as a potential therapeutic strategy \citep{2015_Nassal}.  
		\subsection{Impacts of intracellular delay} \label{Impacts of intracellular delay}
		Time delay is seen to be  an important factor in the intracellular replication process for some viruses. Compared to an ODE model, a delay differential equation model (DDEM) typically has a much more realistic dynamics. Time delay may be responsible for the loss of stability of a steady-state and the oscillation in dynamics. In general, there are two types of delay: (i) pharmacological and (ii) intracellular. The delay between the ingestion of a drug and its appearance within the cells is known as the pharmacological delay. The time-lapse between the infection of a host cell and the discharge of viral particles from the infected cell is known as the intracellular delay \citep{2000_Nelson}. In this study, the intracellular delay is incorporated into  every step of viral life cycle to make the process non-instantaneous. The DDEM is given by in system of equations \eqref{delay_model} in Appendix \hyperref[app:A]{A}, and it is solved for different values of delay ($\tau$). As a result, the simulation demonstrates that the intracellular delay has  very little impact on viral dynamics since there is no significant difference among the solutions. In a recent study,  \citet{2015_Murray} also noted the negligible influence of the delay on the infection dynamics.
		\subsection{Impacts of  double-stranded linear DNA-containing capsids}\label{the effects of double stranded linear DNA-containing capsids}
		The dslDNA-containing capsids are one of the intracellular components of the virus. These can be  integrated in the host cell genome. dslDNA-containing capsids could be released as  virions that contain dslDNAs or could be transported to the nucleus. dslDNAs are  defective form of the viral DNA. These can produce viral surface proteins, but are unable to produce functional pgRNAs due to some mutations that are introduced when These are converted into cccDNAs \citep{2021_tu_hepatitis}.
		   The impacts of  dslDNA-containing capsids on infection dynamics are examined in this study.
		In Section \ref{model_formulation}, the roles of dslDNAs to the production of  2.4 and 2.1 kb mRNAs  are incorporated  into the model \eqref{full model} through the term $\lambda_{sdl} D_L$. Keeping other parameters fixed, the  system of equations \eqref{full model} is solved for different values of $\lambda_{sdl}$ and the solutions are compared. No substantial change  is observed in the compartments of the proposed model especially in the cccDNA and virus compartments. Therefore, targeting the dslDNA containing capsids for possible future treatment options does not seem promising. Moreover, the differential equation corresponding to the dslDNA-containing capsids may be ignored in this model \eqref{full model} for further analysis.
		\section{Global sensitivity analysis of model parameters} \label{GSA}
		Because of uncertainties in the experimental data used in the estimation of model parameters, the accuracy of outputs of a mathematical model, related to specific biological phenomena becomes frequently poor. Most of the  studies generally  focus on local sensitivity analysis (LSA). However, LSA does not provide the complete necessary information about the  uncertainty and sensitivity of the concerned model parameters. In this case, global sensitivity analysis (GSA) performs well and can clearly describe the contributions of each model parameter  irrespective of the roles of other parameters. The GSA is a statistical technique, which is used to study the sensitivity of parameters of a system or of a mathematical model. Various methods are  used for the global sensitivity analysis, such as Sobol indices, Fourier amplitude sensitivity test, partial rank correlation coefficient (PRCC). In this study, Latin hypercube sampling-partial rank correlation coefficient (LHS-PRCC) method is applied. This method is well-explained in the article of  \citet{2008_marino_methodology}.   In this method, PRCC values can provide relevant useful information. PRCC can also aid us in determining the most influential set of parameters to achieve specific objectives in the elimination of this disease. 
		\subsection{Scatter plots: The monotonic relationship between input and output variables} \label{Result: Scatter plots}
\begin{figure}[h!]
		\centering
		{(A): Scatter plot for rcDNA ($R$)}
		\includegraphics[height=3.4cm,width=16cm]{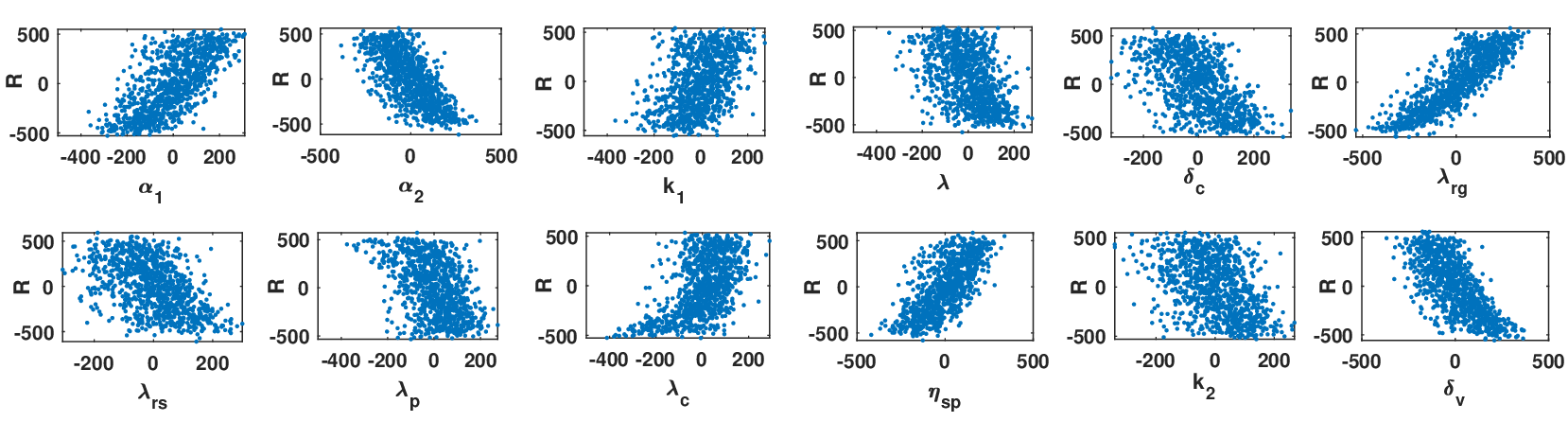}			
		\par\noindent\rule{\textwidth}{1pt}
		{(B): Scatter plots for cccDNA ($C$)}
		\includegraphics[height=3.4cm,width=16cm]{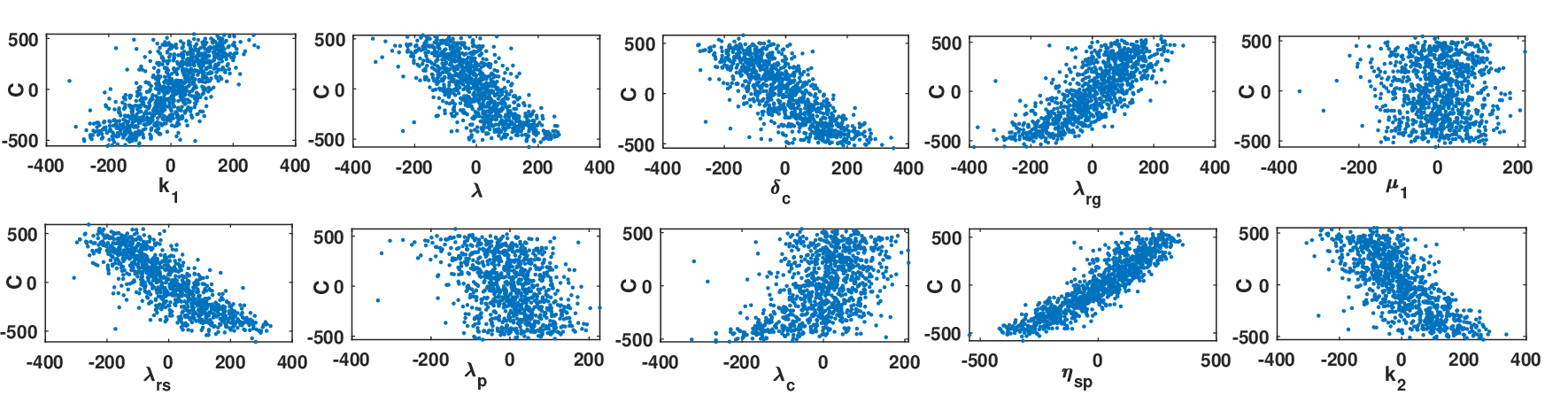}
		\par\noindent\rule{\textwidth}{1pt}
		{(C): Scatter plots for 3.5 kb pgRNA ($R_g$)}
		\includegraphics[height=3.4cm,width=16cm]{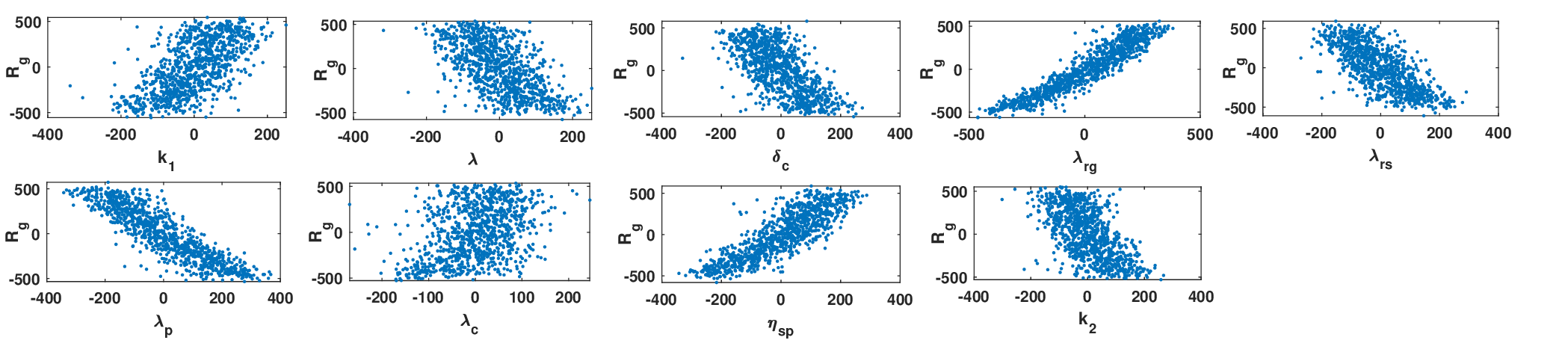}
		\par\noindent\rule{\textwidth}{1pt}	
		{(D): Scatter plots for 2.4 and 2.1 kb mRNA ($R_s$)}	
		\includegraphics[height=3.4cm,width=16cm]{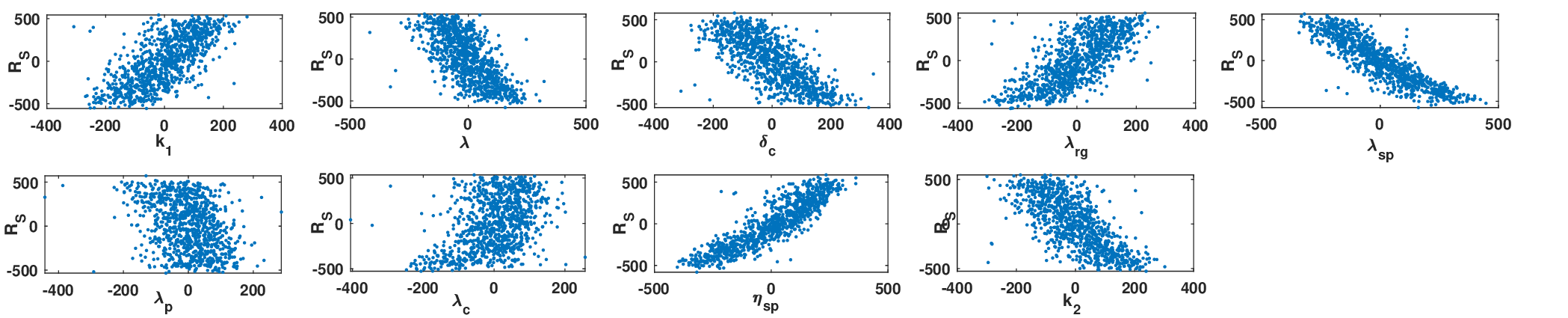}
		\par\noindent\rule{\textwidth}{1pt}	
		{(E): Scatter plots for 0.7 kb mRNA ($R_h$)}\\	
		\includegraphics[height=3.1cm,width=16cm]{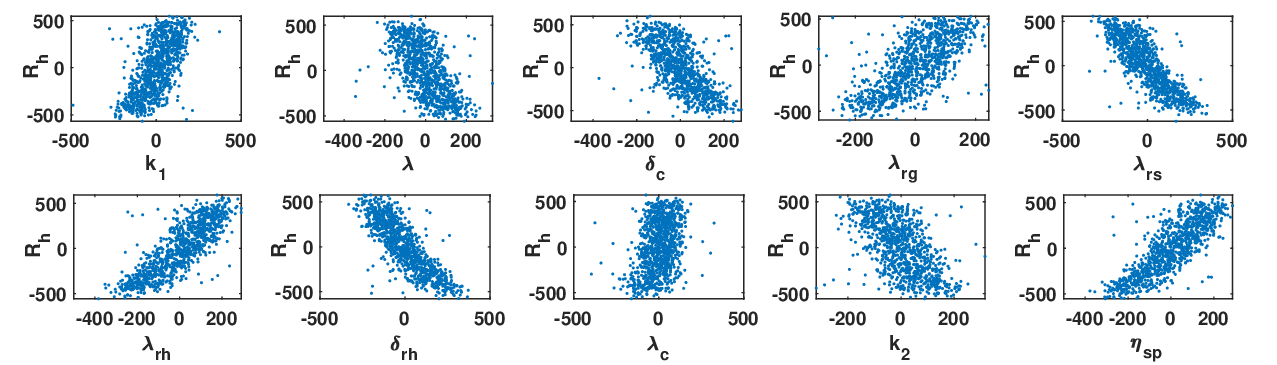}	
		\par\noindent\rule{\textwidth}{1pt}	
		{(F): Scatter plots for HBx protein ($H$)}
			\includegraphics[height=3.1cm,width=16cm]{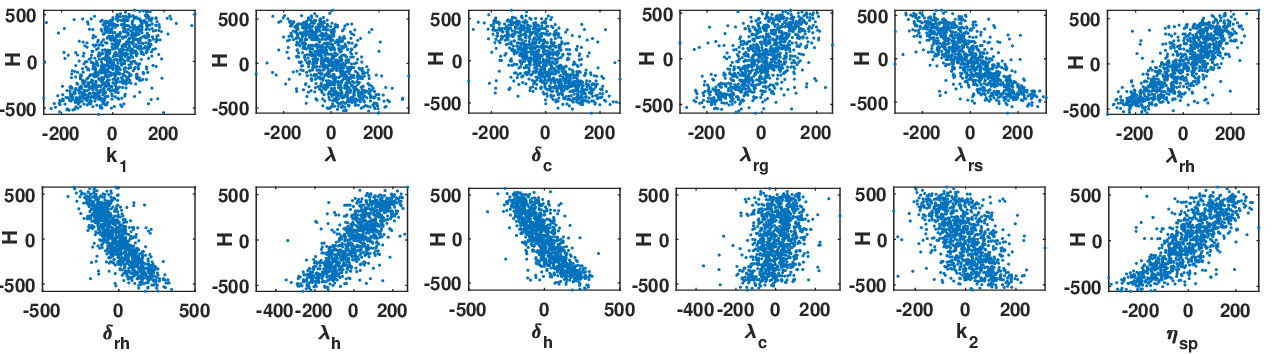}	
		\caption{Scatter plots of rcDNA, cccDNA, 2.4 and 2.1 kb mRNA, 0.7 kb mRNA, and HBx protein.}
		\label{Scatter plot:  R to Z}
	\end{figure}

\begin{figure}[h!]
		\centering
		{(A): Scatter plots for polymerase ($P$)}\\
		\includegraphics[height=2.0cm,width=8cm]{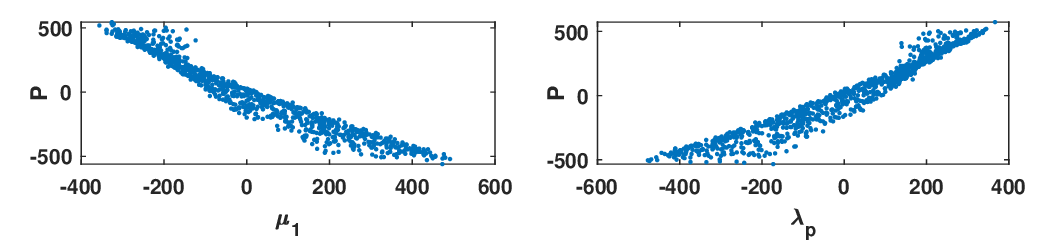}
		\par\noindent\rule{\textwidth}{1pt}	
		{(B): Scatter plots for core protein ($C_p$)}
		\includegraphics[height=1.7cm,width=16cm]{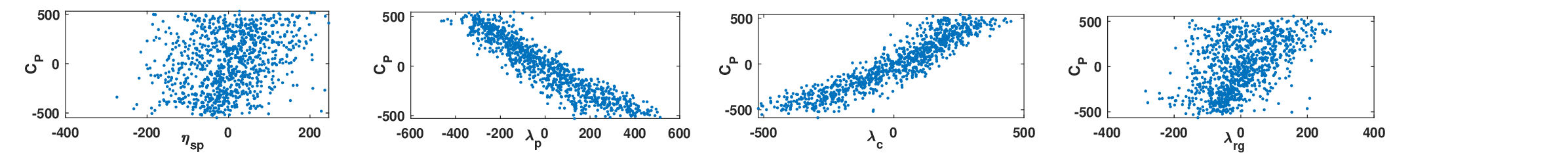}
		\par\noindent\rule{\textwidth}{1pt}	
		{(C): Scatter plots for RNP complex ($Z$)}
			\includegraphics[height=2.1cm,width=11cm]{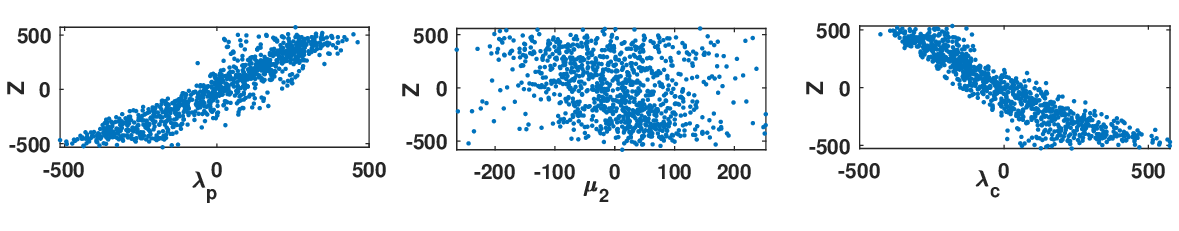}
		\par\noindent\rule{\textwidth}{1pt}
		{(D): Scatter plots for pgRNA-containing capsid ($P_g$)}
		\includegraphics[height=3.5cm,width=16cm]{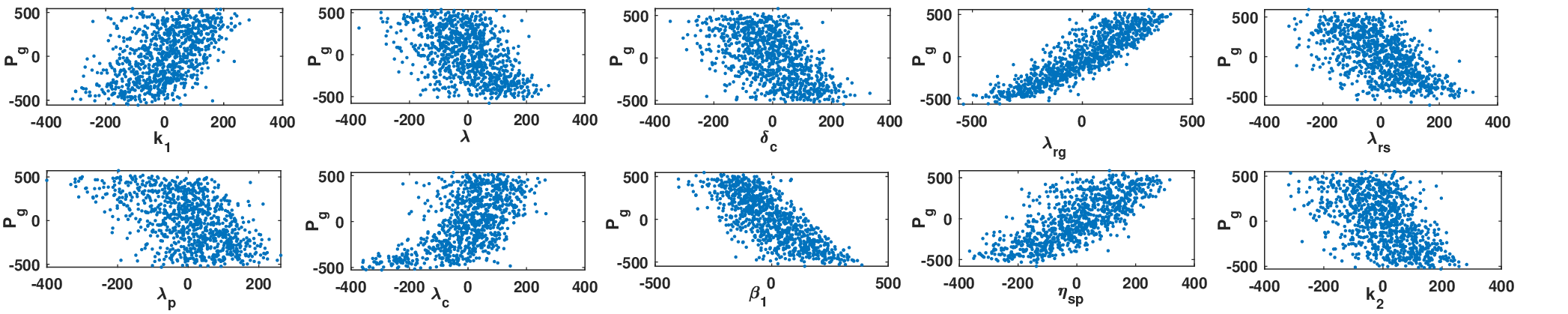}
		\par\noindent\rule{\textwidth}{1pt}
		{(E): Scatter plots for surface protein ($S_p$)}
		\includegraphics[height=3.5cm,width=16cm]{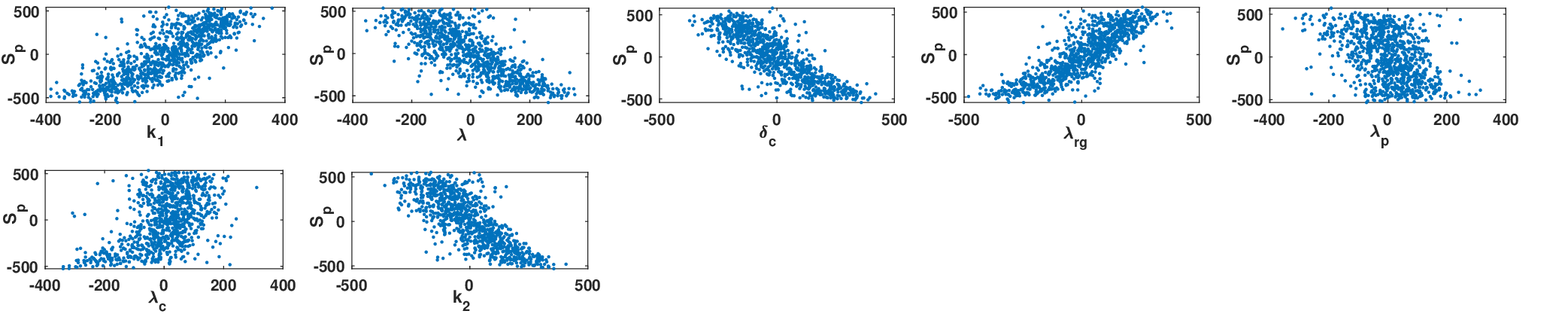}
		\par\noindent\rule{\textwidth}{1pt}
		{(F): Scatter plots for ssDNA-containing capsid ($S$)}
		\includegraphics[height=3.5cm,width=16cm]{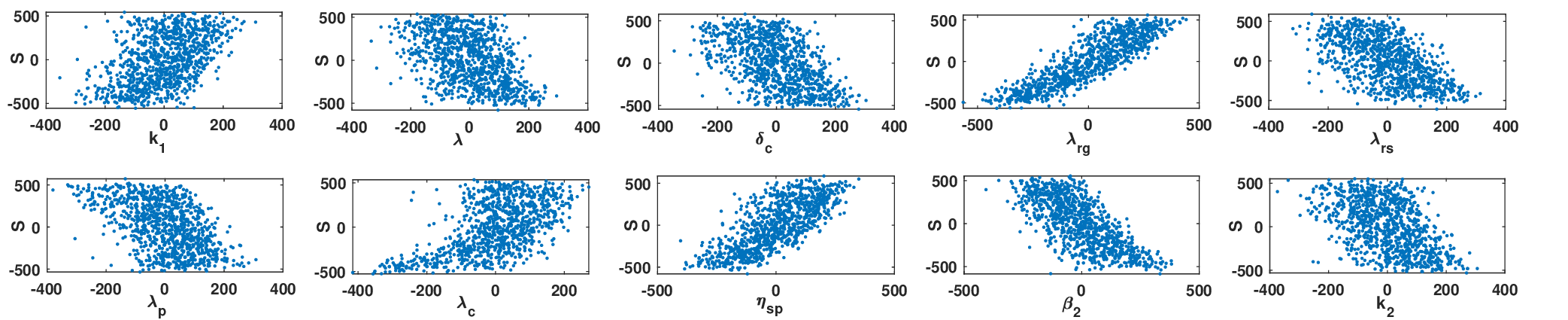}
		\par\noindent\rule{\textwidth}{1pt}
		\caption{Scatter plots of polymerase, core protein, RNP complex,  pgRNA-containing capsid, surface protein, and ssDNA-containing capsid.}
		\label{Scatter plot:  Pg to V}
	\end{figure}

\begin{figure}[h!]
		\centering
		{(A): Scatter plots for dslDNA-containing capsid ($D_L$)}
		\includegraphics[height=3.8cm,width=16cm]{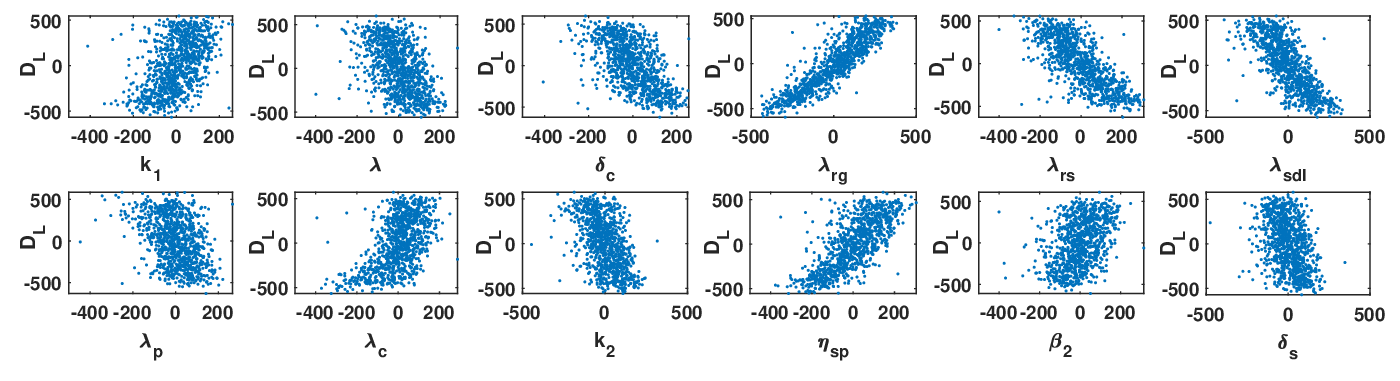}\
		\par\noindent\rule{\textwidth}{1pt}
		{(B): Scatter plots for dsDNA-containing capsid ($D$)}
		\includegraphics[height=3.4cm,width=16cm]{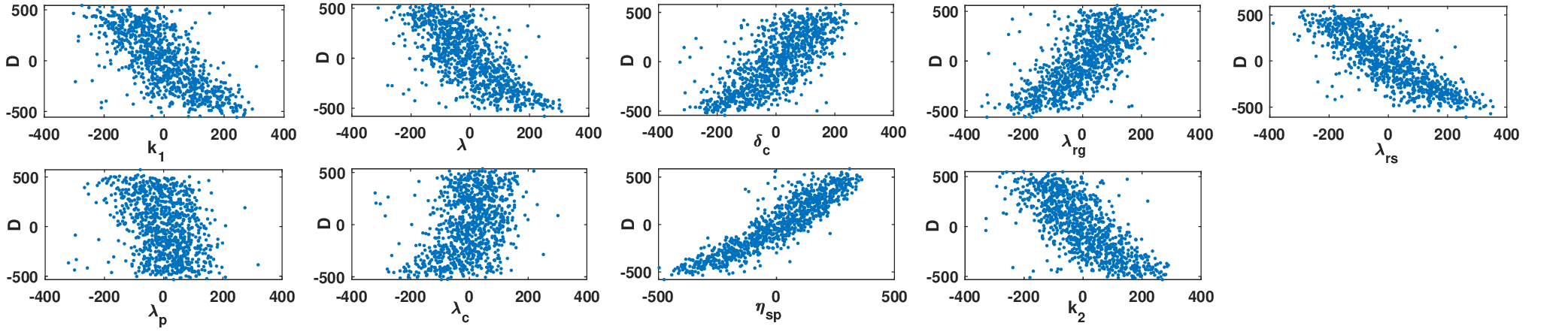}\
		\par\noindent\rule{\textwidth}{1pt}
		{(C): Scatter plots for viruse ($V$)}
		\includegraphics[height=3.7cm,width=16cm]{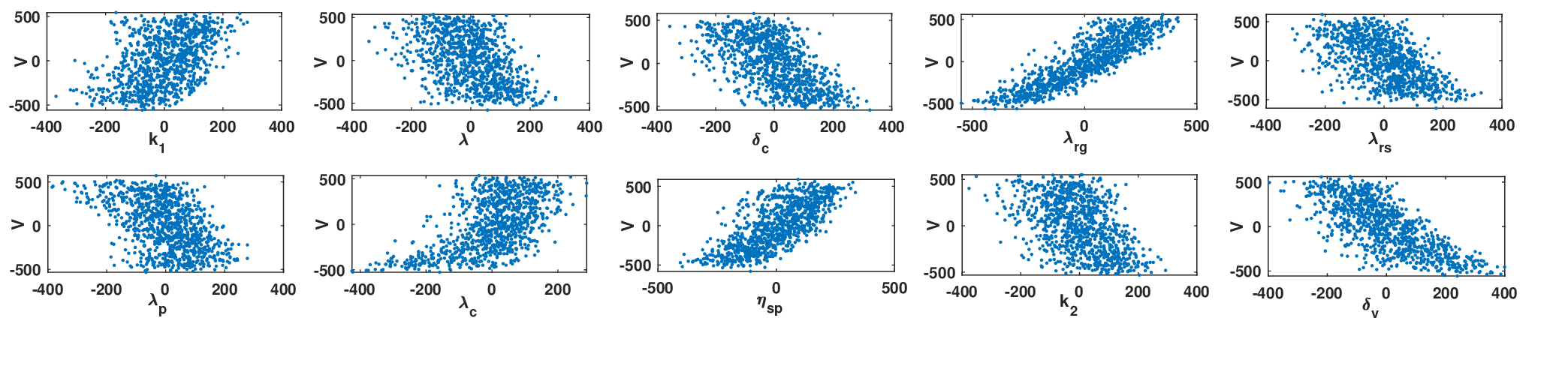}
			\par\noindent\rule{\textwidth}{1pt}
		{(D): Scatter plots for SVP ($S_v$)}
		\includegraphics[height=3.7cm,width=16cm]{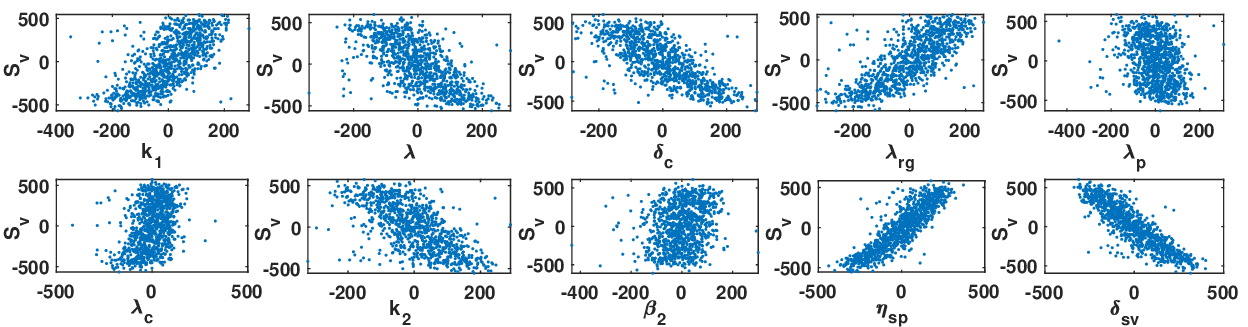}
		\caption{Scatter plots of dslDNA-containing capsid, dsDNA-containing capsid, virus and SVP.}
		\label{Scatter plot:  DL to SVPs}
	\end{figure}

		In order to make better prediction on the infection dynamics and to recommend  new different type of treatment strategies, it is  essential to explore  how the outputs of a system are influenced if the  values of the associated parameters vary within  biologically relevant ranges.
		%
		Simulation results of the full  model \eqref{full model} are shown by scatter plots which are presented in Figures \ref{Scatter plot:  R to Z}, \ref{Scatter plot:  Pg to V} and \ref{Scatter plot:  DL to SVPs}.   
		PRCC values of all model  parameters are computed  with respect to  the dependent variables  on the
		 90$^{th}$ day and are visualized in  Table \ref{PRCC-table}.
		The positive PRCC value of a model parameter and a compartment of the model indicates that any increase or decrease in the parameter's value, whether individually or simultaneously, leads to an  enhancement or reduction  in the concentration of that compartment.  
		On the other hand, negative correlation (PRCC value negative) shows the opposite aspects.
		Based on the PRCC values, the most positively sensitive (MPS), the most negatively sensitive (MNS) and  less sensitive parameters are identified and listed  in   Table \ref{Complete list of positively correlated, negatively correlated}.
		
		Global sensitivity analysis reveals several  striking results. Based on the PRCC values (Table \ref{PRCC-table}) and the outputs shown in Table \ref{Complete list of positively correlated, negatively correlated}, some of the findings are listed below.
		\begin{enumerate}[(i)]
			\item  The production rate of SVPs ($\eta_{sp}$) is positively correlated with almost all viral components mentioned in the model \eqref{full model} except surface proteins.
		It means that the production rate  of SVPs  may play a key role in enhancing the infection. In this context, HBsAg release inhibitors, such as REP 2055 and REP 2139-Ca \citep{2016_AI-Mahtab}, and REP 2139-Mg and REP 2165-Mg \citep{2020_Bazinet} (both are in clinical trial), may be potentially useful components of future therapies. Besides, REP 9AC is also designed as a potent HBsAg release inhibitor \citep{2011_Mahtab}.
		
			\item  All the  components except dsDNA-containing capsids have  positive correlation with  recycling rate $(k_1)$. Through the recycling of capsids, newly produced capsids  are transported  to the nucleus  from the cytoplasm, thereby increasing the pool of cccDNAs. 
		Compared to other components, the surface protein ($S_p$) is the most sensitive components for this parameter, exhibiting  PRCC value of 0.8067. Therefore,  it can be considered that capsid recycling may enhance the infection and serve as a positive feedback loop. 
			 Long-term treatment with nucleos(t)ide analogs (NAs) can decrease the cccDNA pool in HBV-infected hepatocytes by inhibiting nucleocapsid recycling \citep{2022_Leowattana}. Besides, \citet{2022_wong} observed that capsid assembly
			 modulators (CAMs) can inhibit cccDNA replenishment by disrupting capsid disassembly and nucleocapsid recycling. So, targeting capsid recycling could be  an effective strategy to eradicate the infection.
			
			\item The production rate of  3.5 kb pgRNAs ($\lambda_{rg}$)   is positively associated with all viral compartments considered in the model \eqref{full model}, implying that  transcription of cccDNAs is a critical step in the viral life cycle. 	
			 On the other hand, transcription rate of 2.4 kb and 2.1 kb mRNAs $(\lambda_{rs})$ is negatively correlated with nearly all compartments except 2.4 kb and 2.1 kb mRNA, surface protein and virus. It indicates that the availability of surface proteins leads to a higher quantity of viruses, resulting enhancement of the infection. This accurately reflects the biological mechanisms of the viral infection. Hence, transcription of cccDNAs should be considered while proposing any new control strategy. Disruption of cccDNAs is one of the ways to reduce the infection. In Subsection \ref{Effects of initial concentration  of cccDNA}, it is seen that the presence of a
			 few copies of cccDNAs  can re-initiate and blow-up the infection.  Gene editing drugs like clustered regularly interspaced short palindromic repeats-associated 9 (CRISPR/Cas9 \citep{2021_Stone}) and zinc-finger nucleases  (ZFNs \citep{2014_Weber}) (both are in clinical trial) may show strong potential as standalone therapy or as synergy therapy. Apart from this, in order to reduce HBV transcription,  transcription inhibitors (ARC-520 \citep{2020_Yuen}, GSK3389404 \citep{2019_Han}) and gene expression inhibitor (RG7834 \citep{2020_Menne}) can be considered in the future.
			\item It is also observed that the production rate of polymerases $(\lambda_p)$  and the production rate of core proteins $(\lambda_{c})$  both are  sensitive (positive for some compartments and negative for others) for all compartment of the model. In addition, these two parameters are synchronously MPS and MNS. This kind of behavior is quite peculiar.  As far as the infection is concerned,  production rates of polymerases and core proteins play dual role (simultaneously MPS and MNS). This is one of the crucial findings in this study. Probably, dual role of these two parameters are noticed in this study for the first time. No one has informed about this kind of behaviors of these two parameters so far. In order to target capsid assembly, it is essential to disrupt the formation of capsid-like structures. In this scenario, the capsid assembly modulator drug NVR3-778 \citep{2020_Yuen_Safety} (under clinical trials) can be focused. 
			
			\item  Since the production rates of ssDNA and dsDNA-containing capsids ($\beta_1$ and $\beta_2$) promote infection, nucleos(t)ide analogs (NAs) such as lamivudine (LAM), telbivudine (LDT), entecavir (ETV), adefovir dipivoxil (ADV), and tenofovir disoproxil fumarate (TDF) \citep{2022_Leowattana}  are promising candidates for inhibiting the infection. Nowadays, these drugs are widely used due to their effectiveness.

			\item Global sensitivity analysis identifies that the production rate of 3.5 kb pgRNA ($\lambda_{rg}$) plays as the most positively influential parameter in the framework of  the system delineated by equation \eqref{full model}. This is because  it is  the most positively sensitive parameter, as reflected by its high PRCC values, among the following compartments: rcDNA-containing capsids, 3.5 kb pgRNA, pgRNA-containing capsids, surface proteins, ssDNA-containing capsids, and viruses.
		\end{enumerate}
		%
		It is important to note that  all the drugs which are mentioned above except NAs  are in clinical trial and may be considered for use in the future as monotherapy or combination therapy (CT), subject to approval by the Food and Drug Administration (FDA). In this study, we do not account the functions of  interferon therapy (IFN-$\alpha$ or Peg-IFN-$\alpha$) because the proposed model does not include immune cells. Considering all these results and findings, it is believed that  the best way to control this disease as well as the best CT  can be determined.

        In order to reduce the computational cost and make the model easier to handle both analytically
and computationally, a simplified version of the proposed model \eqref{full model}  is given in Appendix \hyperref[app:B]{B}.
		\begin{landscape}
			\begin{table}[]
				\small
				\begin{tabular}{|l|c|c|c|c|c|c|c|c|c|c|c|c|c|c|c|c|}
					\hline
					\diagbox[width=1.1cm]{P}{V}            & $R$                             & $C$                             & $R_g$                           & $R_s$                           & $R_h$                           & $H$                             & $P$                             & $Z$                             & $C_p$                           & $P_g$                           & $S_p$                           & $S$                             & $D_L$                           & $D$                             & $V$                             & $S_v$                           \\ \hline
					$\alpha_1$      & 0.7028                          & 0.0765                          & 0.0702                          & 0.0652                          & 0.0779                          & 0.0257                          & -0.0074                         & -0.0496                         & 0.0709                          & 0.0626                          & 0.0502                          & 0.0409                          & 0.0446                          & -0.0909                         & 0.0529                          & 0.0899                          \\ \hline
					$\alpha_2$      & -0.5793                         & 0.0556                          & 0.0184                          & 0.0158                          & 0.0434                          & -0.0327                         & -0.0424                         & -0.0254                         & 0.0459                          & 0.0232                          & 0.0050                          & 0.0146                          & -0.0020                          & 0.0013                          & 0.0048                          & -0.0279                         \\ \hline
					$\delta_r$      & -0.2086                         & -0.0080                          & 0.0045                          & 0.0141                          & -0.0122                         & -0.0097                         & 0.0119                          & -0.0278                         & 0.0201                          & -0.0080                          & -0.0315                         & -0.0211                         & 0.0129                          & 0.0136                          & -0.0089                         & -0.0218                         \\ \hline
					$k_1$           & 0.488                           & 0.7193                          & 0.6203                          & 0.7285                          & \cellcolor[HTML]{FFFFFF}0.6489  & 0.6000                          & 0.1035                          & 0.1941                          & 0.1025                          & 0.5278                          & 0.7847                          & 0.5439                          & 0.5217                          & -0.6940                          & 0.5175                          & 0.6815                          \\ \hline
					$\lambda$       & -0.4193                         & -0.6984                         & -0.6096                         & -0.6917                         & -0.6307                         & -0.5692                         & -0.0955                         & -0.1773                         & -0.1195                         & -0.5290                          & -0.7758                         & -0.5320                          & -0.5150                          & -0.7192                         & -0.4540                          & -0.6739                         \\ \hline
					$\delta_c$      & -0.6083                         & -0.7720                          & -0.6918                         & -0.7735                         & -0.7083                         & -0.6672                         & -0.1768                         & -0.2298                         & -0.1121                         & -0.6042                         & \cellcolor[HTML]{FE0000}-0.8383 & -0.6253                         & -0.6009                         & 0.6315                          & -0.6574                         & -0.7477                         \\ \hline
					$\lambda_{rg}$  & \cellcolor[HTML]{34CDF9}0.8681  & 0.7812                          & \cellcolor[HTML]{34CDF9}0.9307  & 0.7757                          & 0.7292                          & 0.6843                          & 0.4414                          & 0.4742                          & 0.3955                          & \cellcolor[HTML]{34CDF9}0.8965  & \cellcolor[HTML]{34CDF9}0.8420  & \cellcolor[HTML]{34CDF9}0.9055  & \cellcolor[HTML]{34CDF9}0.8932  & 0.7588                          & \cellcolor[HTML]{34CDF9}0.8888  & 0.7675                          \\ \hline
					$\mu_1$         & -0.0587                         & -0.0717                         & -0.1929                         & -0.0686                         & -0.0661                         & -0.0171                         & \cellcolor[HTML]{FE0000}-0.9570 & 0.1711                          & -0.1663                         & -0.0709                         & -0.0603                         & -0.0799                         & -0.0529                         & -0.0269                         & -0.0756                         & -0.0456                         \\ \hline
					$\delta_{rg}$   & -0.0076                         & -0.0183                         & 0.0019                          & 0.0024                          & -0.0254                         & -0.0364                         & 0.0130                          & 0.0596                          & -0.0865                         & -0.0270                          & -0.0089                         & 0.0006                          & -0.0375                         & -0.0110                          & -0.0039                         & -0.0199                         \\ \hline
					$\lambda_{rs}$  & -0.6813                         & \cellcolor[HTML]{FE0000}-0.8897 & -0.8317                         & 0.0497                          & -0.8485                         & -0.8157                         & -0.3100                         & -0.2260                          & -0.3142                         & -0.7689                         & -0.0489                         & -0.7906                         & -0.7634                         & \cellcolor[HTML]{FE0000}-0.8969 & -0.7105                         & 0.0535                          \\ \hline
					$\lambda_{sdl}$ & 0.0034                          & 0.0032                          & 0.0157                          & -0.0175                         & 0.0171                          & -0.0016                         & 0.0029                          & -0.0061                         & 0.0318                          & 0.0328                          & -0.0080                          & 0.0475                          & \cellcolor[HTML]{FE0000}-0.7637 & 0.0133                          & 0.0151                          & -0.0265                         \\ \hline
					$\lambda_{sp}$  & 0.0054                          & 0.0320                           & -0.0140                          & \cellcolor[HTML]{FE0000}-0.8956 & 0.0051                          & -0.0029                         & 0.0207                          & -0.0365                         & 0.0355                          & 0.0060                          & -0.0173                         & -0.0189                         & 0.0075                          & 0.0915                          & -0.0121                         & -0.0160                          \\ \hline
					$\delta_{rs}$   & 0.0289                          & 0.0319                          & 0.0050                           & 0.0409                          & 0.0632                          & 0.0792                          & 0.0541                          & 0.0439                          & -0.0369                         & 0.0091                          & -0.0002                         & 0.0400                            & 0.0102                          & -0.0013                         & 0.0153                          & 0.0284                          \\ \hline
					$\lambda_{rh}$  & 0.0066                          & 0.0120                           & -0.0071                         & 0.0757                          & \cellcolor[HTML]{34CDF9}0.8539  & 0.8201                          & 0.0581                          & 0.0140                           & -0.0172                         & -0.0124                         & 0.0208                          & 0.0035                          & -0.0279                         & -0.0324                         & 0.0008                          & 0.0385                          \\ \hline
					$\delta_{rh}$   & -0.0019                         & 0.0232                          & 0.0038                          & 0.0400                            & \cellcolor[HTML]{FE0000}-0.8594 & \cellcolor[HTML]{FE0000}-0.8297 & 0.0020                           & -0.0649                         & 0.0606                          & 0.0142                          & -0.0265                         & -0.0274                         & -0.0293                         & 0.0293                          & 0.0021                          & 0.0176                          \\ \hline
					$\lambda_h$     & 0.0080                          & 0.0171                          & 0.0064                          & 0.0131                          & 0.0338                          & \cellcolor[HTML]{34CDF9}0.8337  & 0.0033                          & -0.0217                         & 0.0353                          & 0.0099                          & 0.0501                          & 0.0172                          & 0.0116                          & 0.0168                          & 0.0006                          & -0.0032                         \\ \hline
					$\delta_h$      & -0.0153                         & -0.0396                         & -0.0365                         & -0.0188                         & -0.0046                         & -0.8293                         & -0.0752                         & -0.0052                         & 0.0051                          & -0.0297                         & -0.0309                         & -0.0098                         & -0.0312                         & 0.0049                          & -0.0375                         & 0.0152                          \\ \hline
					$\lambda_p$     & -0.3890                         & -0.2954                         & \cellcolor[HTML]{FE0000}-0.8756 & -0.3365                         & -0.2385                         & -0.2550                         & \cellcolor[HTML]{34CDF9}0.9561  & \cellcolor[HTML]{34CDF9}0.8981  & \cellcolor[HTML]{FE0000}-0.8862 & -0.4508                         & -0.3448                         & -0.4607                         & -0.4302                         & -0.3000                         & -0.4111                         & -0.2972                         \\ \hline
					$\delta_p$      & 0.0322                          & 0.0612                          & 0.1754                          & 0.1076                          & 0.0324                          & 0.0534                          & -0.2745                         & -0.1825                         & 0.1653                          & 0.0341                          & 0.0362                          & 0.0186                          & 0.0433                          & 0.0271                          & 0.0316                          & 0.0722                          \\ \hline
					$\mu_2$         & -0.013                          & -0.0394                         & -0.0395                         & -0.0047                         & -0.0214                         & 0.0199                          & -0.0076                         & -0.4703                         & -0.1125                         & -0.0135                         & 0.0146                          & -0.0145                         & -0.0154                         & -0.0127                         & -0.0139                         & -0.0056                         \\ \hline
					$\delta_z$      & -0.0099                         & -0.0221                         & 0.0050                           & -0.0458                         & -0.0377                         & -0.0519                         & -0.1105                         & -0.0820                          & 0.0492                          & 0.0023                          & -0.0022                         & -0.0081                         & -0.0184                         & -0.0189                         & -0.0155                         & -0.0161                         \\ \hline
					$\lambda_c$     & 0.5121                          & 0.4449                          & 0.3412                          & 0.4400                            & 0.3268                          & 0.3366                          & 0.0929                          & \cellcolor[HTML]{FE0000}-0.9147 & \cellcolor[HTML]{34CDF9}0.9074  & 0.5736                          & 0.4753                          & 0.5956                          & 0.5683                          & 0.4035                          & 0.5525                          & 0.4132                          \\ \hline
					$\delta_{cp}$   & -0.0443                         & -0.0629                         & -0.0425                         & -0.0296                         & -0.0379                         & -0.0464                         & -0.0113                         & 0.4331                          & -0.3823                         & -0.0515                         & -0.0653                         & -0.0618                         & -0.0416                         & -0.0469                         & -0.0440                          & -0.0284                         \\ \hline
					$\beta_1$       & 0.0747                          & 0.0612                          & 0.0585                          & 0.0435                          & 0.0344                          & 0.0167                          & -0.0528                         & 0.0076                          & 0.0140                          & \cellcolor[HTML]{FE0000}-0.7660 & 0.0322                          & 0.0665                          & 0.0711                          & 0.0354                          & 0.0899                          & 0.0864                          \\ \hline
					$\delta_{pg}$   & -0.0152                         & -0.0289                         & -0.0307                         & -0.0156                         & -0.0394                         & -0.0714                         & -0.0236                         & -0.0384                         & 0.0187                          & -0.0062                         & -0.0393                         & -0.0315                         & -0.0261                         & -0.0139                         & -0.0174                         & -0.0158                         \\ \hline
					$k_2$           & -0.3887                         & -0.6931                         & -0.6097                         & -0.6969                         & -0.6412                         & -0.5926                         & -0.1360                         & -0.1860                         & -0.1088                         & -0.5234                         & -0.7830                         & -0.5360                         & -0.5063                         & -0.7083                         & -0.4464                         & -0.6838                         \\ \hline
					$\eta_{sp}$     & 0.6743                          & \cellcolor[HTML]{34CDF9}0.8932  & 0.8323                          & \cellcolor[HTML]{34CDF9}0.8989  & 0.8545                          & 0.8300                          & 0.3079                          & 0.3114                          & 0.2530                          & 0.7748                          & -0.0048                         & 0.7945                          & 0.7784                          & \cellcolor[HTML]{34CDF9}0.8977  & 0.7149                          & \cellcolor[HTML]{34CDF9}0.8850  \\ \hline
					$\delta_{sp}$   & 0.0298                          & 0.0626                          & 0.0253                          & 0.0813                          & 0.0308                          & 0.0362                          & 0.0392                          & 0.0195                          & -0.0126                         & 0.0509                          & -0.0170                          & 0.0519                          & 0.0095                          & 0.0900                            & 0.0183                          & 0.0426                          \\ \hline
					$\beta_2$       & 0.3564                          & 0.2536                          & 0.2038                          & 0.2534                          & 0.2256                          & 0.2058                          & 0.1041                          & 0.0759                          & 0.0092                          & 0.1734                          & 0.3518                          & \cellcolor[HTML]{FE0000}-0.6284 & 0.4038                          & 0.2331                          & 0.3821                          & 0.2743                          \\ \hline
					$\delta_s$      & -0.3336                         & -0.2236                         & -0.1634                         & -0.2011                         & -0.2067                         & -0.1687                         & -0.0432                         & -0.0340                          & -0.0277                         & -0.1099                         & -0.2813                         & -0.3902                         & -0.3595                         & -0.2233                         & -0.3669                         & -0.2148                         \\ \hline
					$\delta_{dl}$   & 0.0368                          & -0.004                          & 0.0025                          & 0.0300                          & 0.0284                          & 0.0488                          & 0.0463                          & 0.0354                          & 0.0053                          & 0.0120                          & 0.0088                          & 0.0305                          & 0.0191                          & 0.0332                          & 0.0459                          & -0.0242                         \\ \hline
					$\delta_d$      & -0.1099                         & -0.0972                         & -0.1232                         & -0.1228                         & -0.0662                         & -0.0758                         & 0.0090                          & -0.0020                         & -0.0210                          & -0.1274                         & -0.1347                         & -0.0939                         & -0.1168                         & -0.0887                         & -0.1442                         & -0.1289                         \\ \hline
					$\delta_v$      & \cellcolor[HTML]{FE0000}-0.6902 & -0.0578                         & -0.0113                         & -0.1005                         & -0.0027                         & -0.0024                         & 0.0226                          & -0.0287                         & 0.0203                          & 0.0051                          & -0.0500                         & -0.0371                         & -0.0105                         & 0.1126                          & \cellcolor[HTML]{FE0000}-0.7232 & -0.0504                         \\ \hline
					$\delta_{sv}$   & 0.0012                          & 0.0699                          & 0.0398                          & 0.0242                          & 0.0486                          & 0.0597                          & 0.0066                          & 0.0022                          & 0.0124                          & 0.0523                          & -0.0165                         & 0.0200                            & 0.0471                          & 0.0463                          & 0.0400                          & \cellcolor[HTML]{FE0000}-0.8778 \\ \hline
				\end{tabular}
				\caption{List of PRCC values. As shown in the north-west box, P stands for parameters and V stands for variables. Blue  and red cells represent highest and lowest PRCC values of the corresponding  variable of the model.}
				\label{PRCC-table}
			\end{table}
		\end{landscape}
				\begin{table}[]
			\small
			\caption{A comprehensive categorization of parameters is provided, classified as positively correlated, negatively correlated, less sensitive, the most positively sensitive, and the most negatively sensitive. The positively correlated parameters are listed in descending order, while the negatively correlated parameters are listed in ascending order. Parameters with PRCC values exceeding $0.2$ are identified as positively sensitive, those with PRCC values below $-0.2$ are classified as negatively sensitive, and parameters with PRCC values within the range $\left[-0.2, 0.2\right]$ are deemed less sensitive.
			}
			\begin{tabular}{|c|c|c|c|c|c|}
				\hline
				Variable & \begin{tabular}[c]{@{}c@{}}Positively  correlated \\ parameters\\ PRCC $>0.2$\end{tabular}                             & \begin{tabular}[c]{@{}c@{}}Negatively correlated \\ parameters\\ PRCC$<-0.2$\end{tabular}                                            & \begin{tabular}[c]{@{}c@{}}Less\\ sensitive \\ parameters\\ $-0.2\leq$PRCC$\leq 0.2$\end{tabular} & \begin{tabular}[c]{@{}c@{}}Most \\ positively \\ sensitive \\ paramter\end{tabular} & \begin{tabular}[c]{@{}c@{}}Most\\  negatively\\   sensitive\\  paramater\end{tabular} \\ \hline
				$R$      & \begin{tabular}[c]{@{}c@{}}$\lambda_{rg}, \alpha_1, \eta_{sp}, \lambda_c,$\\  $k_1, \beta_2$\end{tabular}                & \begin{tabular}[c]{@{}c@{}}$\delta_v, \lambda_{rs}, \delta_c, \alpha_2, \lambda,$\\  $\lambda_p, k_2, \delta_s, \delta_r$\end{tabular} & \begin{tabular}[c]{@{}c@{}}All other \\ parameters\end{tabular}                                   & \cellcolor[HTML]{34CDF9}$\lambda_{rg}$                                              & $\delta_v$                                                                            \\ \hline
				$C$      & $\eta_{sp}, \lambda_{rg}, k_1, \lambda_c, \beta_2$                                                                     & $\lambda_{rs}, \delta_c, \lambda, k_2, \lambda_p, \delta_s$                                                                          & \begin{tabular}[c]{@{}c@{}}All other \\ parameters\end{tabular}                                   & \cellcolor[HTML]{FD6864}$\eta_{sp}$                                                 & $\lambda_{rs}$                                                                        \\ \hline
				$R_{g}$  & $\lambda_{rg}, \eta_{sp}, k_1, \lambda_c, \beta_2$                                                                     & $\lambda_p, \lambda_{rs}, \delta_c, k_2, \lambda$                                                                                    & \begin{tabular}[c]{@{}c@{}}All other \\ parameters\end{tabular}                                   & \cellcolor[HTML]{34CDF9}$\lambda_{rg}$                                              & $\lambda_p$                                                                           \\ \hline
				$R_s$    & $\eta_{sp}, \lambda_{rg}, k_1, \lambda_c, \beta_2$                                                                     & $\lambda_{sp}, \delta_c, k_2, \lambda, \lambda_p, \delta_s$                                                                          & \begin{tabular}[c]{@{}c@{}}All other \\ parameters\end{tabular}                                   & \cellcolor[HTML]{FD6864}$\eta_{sp}$                                                 & $\lambda_{sp}$                                                                        \\ \hline
				$R_h$    & \begin{tabular}[c]{@{}c@{}}$\eta_{sp}, \lambda_{rh}, \lambda_{rg}, k_1,$\\  $\lambda_c, \beta_2$\end{tabular}            & \begin{tabular}[c]{@{}c@{}}$\delta_{rh}, \lambda_{rs}, \delta_c, k_2,$ \\ $\lambda, \lambda_p, \delta_s$\end{tabular}                  & \begin{tabular}[c]{@{}c@{}}All other \\ parameters\end{tabular}                                   & \cellcolor[HTML]{FD6864}$\eta_{sp}$                                                 & $\delta_{rh}$                                                                         \\ \hline
				$H$      & \begin{tabular}[c]{@{}c@{}}$\lambda_h, \eta_{sp}, \lambda_{rh}, \lambda_{rg},$ \\ $k_1, \lambda_c, \beta_2$\end{tabular} & \begin{tabular}[c]{@{}c@{}}$\delta_{rh}, \delta_h, \lambda_{rs}, \delta_c,$ \\ $k_2, \lambda, \lambda_p$\end{tabular}                  & \begin{tabular}[c]{@{}c@{}}All other \\ parameters\end{tabular}                                   & $\lambda_h$                                                                         & $\delta_{rh}$                                                                         \\ \hline
				$P$      & $\lambda_p, \lambda_{rg}, \eta_{sp}$                                                                                   & $\mu_1, \lambda_{rs}, \delta_p$                                                                                                      & \begin{tabular}[c]{@{}c@{}}All other \\ parameters\end{tabular}                                   & $\lambda_p$                                                                         & $\mu_1$                                                                               \\ \hline
				$Z$      & $\lambda_p, \lambda_{rg}, \delta_{cp}, \eta_{sp}$                                                                      & $\lambda_c, \mu_2, \delta_c, \lambda_{rs}$                                                                                           & \begin{tabular}[c]{@{}c@{}}All other \\ parameters\end{tabular}                                   & $\lambda_p$                                                                         & $\lambda_c$                                                                           \\ \hline
				$C_p$    & $\lambda_c, \lambda_{rg}, \eta_{sp}$                                                                                   & $\lambda_p, \delta_{cp}, \lambda_{rs}$                                                                                               & \begin{tabular}[c]{@{}c@{}}All other \\ parameters\end{tabular}                                   & $\lambda_c$                                                                         & $\lambda_p$                                                                           \\ \hline
				$P_g$    & $\lambda_{rg}, \eta_{sp}, \lambda, k_1$                                                                                & $\lambda_{rs}, \beta_1, \delta_c, \lambda, k_2, \lambda_p$                                                                           & \begin{tabular}[c]{@{}c@{}}All other \\ parameters\end{tabular}                                   & \cellcolor[HTML]{34CDF9}$\lambda_{rg}$                                              & $\lambda_{rs}$                                                                        \\ \hline
				$S_p$    & $\lambda_{rg}, k_1, \lambda_c, \beta_2$                                                                                & $\delta_c, k_2, \lambda, \lambda_p, \delta_s$                                                                                        & \begin{tabular}[c]{@{}c@{}}All other \\ parameters\end{tabular}                                   & \cellcolor[HTML]{34CDF9}$\lambda_{rg}$                                              & $\delta_c$                                                                            \\ \hline
				$S$      & $\lambda_{rg}, \eta_{sp}, \lambda_c, k_1$                                                                              & \begin{tabular}[c]{@{}c@{}}$\lambda_{rs}, \beta_2, \delta_c, k_2, \lambda,$\\  $\lambda_p, \delta_s$\end{tabular}                      & \begin{tabular}[c]{@{}c@{}}All other \\ parameters\end{tabular}                                   & \cellcolor[HTML]{34CDF9}$\lambda_{rg}$                                              & $\lambda_{rs}$                                                                        \\ \hline
				$D_L$    & $\lambda_{rg}, \eta_{sp}, \lambda_c, k_1, \beta_2$                                                                     & \begin{tabular}[c]{@{}c@{}}$\lambda_{sdl}, \lambda_{rs}, \delta_c, \lambda,$ \\$ k_2, \lambda_p, \delta_s$\end{tabular}                & \begin{tabular}[c]{@{}c@{}}All other \\ parameters\end{tabular}                                   & \cellcolor[HTML]{34CDF9}$\lambda_{rg}$                                              & $\lambda_{sdl}$                                                                       \\ \hline
				$D$      & $\eta_{sp}, \lambda_{rg}, \delta_c, \lambda, \beta_2$                                                                  & $\lambda_{rs}, \lambda, k_2, k_1, \lambda_p, \delta_s$                                                                               & \begin{tabular}[c]{@{}c@{}}All other \\ parameters\end{tabular}                                   & \cellcolor[HTML]{FD6864}$\eta_{sp}$                                                 & $\lambda_{rs}$                                                                        \\ \hline
				$V$      & $\lambda_{rg}, \eta_{sp}, \lambda, k_1, \beta_2$                                                                       & \begin{tabular}[c]{@{}c@{}}$\delta_v, \lambda_{rs}, \delta_c, \lambda,$ \\ $k_2, \lambda_p, \delta_s$\end{tabular}                     & \begin{tabular}[c]{@{}c@{}}All other \\ parameters\end{tabular}                                   & \cellcolor[HTML]{34CDF9}$\lambda_{rg}$                                              & $\delta_v$                                                                            \\ \hline
				$S_v$    & $\eta_{sp}, \lambda_{rg}, k_1, \lambda_c, \beta_2$                                                                     & $\delta_{sv}, \delta_c, k_2, \lambda, \lambda_p, \delta_s$                                                                           & \begin{tabular}[c]{@{}c@{}}All other \\ parameters\end{tabular}                                   & \cellcolor[HTML]{FD6864}$\eta_{sp}$                                                 & $\delta_{sv}$                                                                         \\ \hline
			\end{tabular}
			\label{Complete list of positively correlated, negatively correlated}
		\end{table}
		\clearpage	
		\section{Conclusions}
		In case of viral infection, intracellular dynamics model has revealed  various intrinsic biological phenomena of an individual infected cell.
		Based on the available  biological and clinical findings about the HBV intracellular dynamics, this study proposes an intracellular dynamics model. 
		The  model solutions are compared with the experimental data obtained from  four humanized mice, and it is seen that  
		the model solutions closely match the experimental data.
		In order to study the global  sensitivity of the model parameters, a sampling-based method, Latin hypercube sampling-partial rank correlation coefficient is applied. The most positively and most negatively sensitive parameters for each compartment of the model, as well as for the entire system, are identified.  It is also uncovered that two parameters ($\lambda_{c},\lambda_{p}$) have  dual role in modulating the  dynamics of the infection since  these two parameters are
		the most positively sensitive for some compartments, at the same time   being the most negatively sensitive for some other  compartments. The findings of this study are expected to provide new insights on  progression of the disease.
		The outcomes of this study lead to the following conclusions:
		\begin{enumerate}[(i)]
			\item  In Section \ref{effect of HBx proteins}, it is
             observed that the numbers of cccDNAs and viruses increase rapidly when HBx is taken into account, \textit{i.e.}, HBx proteins  contribute to the progression of HBV infection.  So, targeting the HBx protein as a  future  antiviral therapy could be a promising strategy to control the infection.
			\item Based on the simulation outcomes shown in  Sections \ref{Impacts of intracellular delay} and \ref{the effects of double stranded linear DNA-containing capsids}, it is realized that the intracellular delay and the dslDNA-containing capsids seem to have little impact on infection dynamics, as the differences between the solution with and without these factors are not substantial. However, further investigation is required to confirm this observation.
			
				
            
            \item The PRCC values of the production rate of SVPs $(\eta_{sp})$ with respect to nearly all compartments of the model (equation \eqref{full model}) are positive (shown in Table \ref{PRCC-table}). This implies that an increase in the value of this parameter enhances the concentration levels of almost all viral components considered in the model (equation \eqref{full model}), with the exception of surface proteins. These results suggest that the production rate of SVPs $(\eta_{sp})$ plays a significant role in   advancement of the infection. In this context, HBsAg release inhibitors, such as REP 2055, REP 2139-Ca, REP 2139-Mg, REP 2165-Mg  may be effective antiviral agents in the future. However, all of these drugs are currently undergoing clinical trials.
            
		\item From the  PRCC Table \ref{PRCC-table}, it is observed that  except dsDNA-containing capsids, other viral components are positively correlated with the recycling rate $(k_1)$ of newly produced capsids. The surface protein ($S_p$) is the most sensitive components to variations in this parameter.
            Surface proteins are the most sensitive to the variations in this parameter. This indicates that the recycling of capsids could enhance  the infection. In this case, the available inhibitors (e.g., capsid assembly modulators) associated with capsids recycling mechanism may be used as additional treatment options  to reduce this infection.
			\item 
            The production rate of 3.5 kb pgRNAs ($\lambda_{rg}$)  is seen  to be the most positively sensitive parameter for many viral components, as shown in Table \ref{Complete list of positively correlated, negatively correlated}. Consequently, this parameter can be considered as one of the most positively sensitive parameters  in shaping the dynamics of the infection. Thus, the disruption of cccDNAs may help in controlling the infection. In this instance, in addition to widely used drugs like nucleos(t)ide analogs and interferon therapy, gene editing drugs (CRISPR/Cas9, ZFNs) and gene expression inhibitor (RG7834) could be the viable candidates  in controlling the infection.
		\end{enumerate}
In this study, we have identified the steps in the viral life cycle that could be targeted by the existing antiviral therapies, as well as those therapies that are currently in clinical trials. However, the question of which monotherapy or combination therapy will be the most effective has not yet been addressed. In order to explore this, the current model would need to be modified by incorporating the effects of antiviral therapy and immune system. In addition, we have counted the number of viruses produced from only a single infected hepatocyte. This model is not able to capture the dynamics of the infection in whole liver. Considering the total number of uninfected hepatocytes, infected hepatocytes, and viruses along with this proposed model, an inter-intracellular multiscale model may be developed in the future. Moreover, in this study, the intracellular processes are assumed to be linear, which is a limitation of this model.  This could be incorporated in future studies to capture a broader range of infection dynamics.
		
		
		
		%
		%
		%
		\subsection*{Supplementary Material}
		Figure S1: The effects of HBx proteins on (i) Upcoming rcDNA-containing capsids, (ii) cccDNAs, (iii) 3.5 kb pgRNAs, (iv) 2.4 and 2.1 kb mRNAs, (v) 0.7 kb pgRNAs, (vi) Polymerases, (vii) RNP complexes, and (viii) Core proteins.\\
		Figure S2: The effects of HBx proteins on (i) pgRNA-containing capsids, (ii) Surface proteins, (iii) ssDNA-containing capsids, (iv) dsDNA-containing capsids, (v) dslDNA-containing capsids, (vi) Viruses, and (vii) SVPs.\\
		Section 2: The derivation of the volume fraction function ($\Phi$) for active cccDNAs.
		\subsection*{Abbreviations}
		HBV: Hepatitis B virus, cccDNA: covalently closed circular DNA, dslDNA: double-stranded linear DNA, HCC: hepatocellular carcinoma, rcDNA: relaxed circular DNA. RNP: Ribonucleoprotein, dsDNA: double-stranded DNA, SVPs: Sub-viral particles, ETV: entecavir, LSA: Local sensitivity analysis, GSA: Global sensitivity analysis. PRCC: Partial rank correlation coefficient, MPS: Most positively sensitive, MNS: Most negatively sensitive, 
		\subsection*{Acknowledgments}
		The first author would like to acknowledge the financial support obtained from the Council of Scientific \& Industrial
		Research (CSIR), Goverment of India under the CSIR Fellowship scheme (File No: 09/731(0171)/2019-EMR-I). The first
		author also thanks for the research facilities received from the Department of Mathematics, Indian Institute of Technology
		Guwahati, India.

			\section*{Appendix A: Intracellular delay model} \label{app:A}
			The delay version of  intracellular dynamics model \eqref{full model} is given by
			\begin{equation}\label{delay_model}
			\left.
			\begin{split} 
							\frac{dR}{dt}&= \alpha_1 V(t-\tau)-\alpha_2 R-\delta_r R,\\
								\frac{dC}{dt}&= \alpha_2 R(t-\tau)+k_1e^{-\lambda S_p(t-\tau)} D(t-\tau)-\delta_c C.\\							
							\frac{dR_g}{dt}&= \lambda_{rg} \Phi C(t-\tau)-\mu_1R_g P(t-\tau)-\delta_{r_g}R_g,\\
							\frac{dR_s}{dt}&= \lambda_{rs} \Phi C(t-\tau)+\lambda_{sdl} D_L(t-\tau)-(\lambda_{s_p}+ \delta_{r_s})R_s,\\	
							\frac{dR_h}{dt}&= \lambda_{rh} \Phi C(t-\tau)-\delta_{r_h}R_h,\\	
							\frac{dH}{dt}&= \lambda_h R_h(t-\tau)-\delta_h H,\\
							\frac{dP}{dt}	&= \lambda_p R_g(t-\tau)-\mu_1R_g(t-\tau)P-\delta_p P,\\
							\frac{dZ}{dt}	&=\mu_1R_g(t-\tau)P(t-\tau)-\mu_2 ZC_p(t-\tau)-\delta_z Z,\\	
							\frac{dC_p}{dt}&=\lambda_c R_g(t-\tau)-\mu_2m Z(t-\tau)C_p-\delta_{c_p}C_p,\\
							\frac{dP_g}{dt}&=\mu_2 Z(t-\tau)C_p(t-\tau)-\beta_1P_g-\delta_{p_g}P_g,\\
							\frac{dS_p}{dt}&=\lambda_{s_p} R_s(t-\tau)
							-k_2\left(1-e^{-\lambda S_p}\right)D(t-\tau)S_p-\eta_{sp}S_p-\delta_{s_p}S_p,\\
							\frac{dS}{dt}&=\beta_1 P_g(t-\tau)- \beta_2 S-\delta_sS,\\
							\frac{dD_L}{dt}	&=0.1 \beta_2 S(t-\tau)-\lambda_{sdl}D_L-\delta_{d_L} D_L,\\	
							\frac{dD}{dt}&=0.9\beta_2 S(t-\tau)-k_1e^{-\lambda S_p(t-\tau)} D-k_2\left(1-e^{-\lambda S_p(t-\tau)}\right)D S_p(t-\tau)-\delta_d D,\\
							\frac{dV}{dt}&= k_2\left(1-e^{-\lambda S_p(t-\tau)}\right)D(t-\tau)S_p(t-\tau)-\delta_v V,\\
							\frac{dS_v}{dt}&=\eta_{sp}S_p(t-\tau)-\delta_{sv}S_v.\\
				\end{split}
			\hspace{0.01cm}
			\right\}
				\end{equation}
				Here, $\tau$ denotes the intracellular delay.

	\section*{Appendix B: Simplification of the proposed model \eqref{full model}} \label{app:B}
    The proposed model \eqref{full model} can be simplified to make it more computationally efficient and easier to analyze by ignoring some less important parameters and components based on some assumptions and findings of this study.
	\begin{enumerate}[(i)]
		\item If the intracellular components are degraded rapidly compared to the recruitment rate, then infection will disappear on its own. Usually, this does not happen. Therefore, degradation rates are considered to be too small  and can be ignored except for cccDNAs and virus \citep{2016_jun_nakabayashi}. These parameters are also found to be  less-sensitive in Section \ref{GSA}.
		\item
		It is assumed that all cccDNAs are active. In this case, the volume fraction of active cccDNAs $\Phi$  is taken to be equal to 1. As a result,  the corresponding equations $\dfrac{dR_h}{dt}= \lambda_{rh} \Phi C-\delta_{r_h}R_h$ and $\dfrac{dH}{dt}= \lambda_h R_h-\delta_h H$ that enhance the number of active cccDNAs are ignored. 
		\item The system of equations \eqref{full model} does not seem to be affected significantly by dslDNA-containing capsids, as observed in Section \ref{the effects of double stranded linear DNA-containing capsids}. Thus, the reaction equation corresponding to dslDNA-containing capsids  of the  dynamics model \eqref{full model} is also ignored.
		\item Since  the other compartments of the  system of  equations \eqref{full model}  don't depend on $\dfrac{dS_v}{dt}=\eta_{sp}S_p-\delta_{sv} S_v,$ (equation for SVPs), it is not considered for further analysis.
	\end{enumerate}
	Therefore, on the basis of these assumptions, the full dynamics model \eqref{full model} is reduced to the following: 
	\begin{equation} 
		\left.
		\begin{split}
			\frac{dR}{dt}	&= \alpha_1 V-\alpha_2 R,\hspace{5.6cm}
			\frac{dC}{dt}	= \alpha_2 R+k_1e^{-\lambda S_p} D-\delta_c C,\\
			\frac{dR_g}{dt}&= \lambda_{rg} \Phi C-\mu_1R_g P, \hspace{4.66cm}
			\frac{dR_s}{dt}= \lambda_{rs} \Phi  C-\lambda_{s_p} R_s,\\
			\frac{dP}{dt}	&= \lambda_p R_g-\mu_1R_gP,\hspace{5.1cm}
			\frac{dZ}{dt}	=\mu_1R_gP-\mu_2 ZC_p,\\
			\frac{dC_p}{dt}&=\lambda_c R_g-\mu_2 m ZC_p ,\hspace{4.6cm}
			\frac{dP_g}{dt}=\mu_2 ZC_p-\beta_1 P_g,\\
			\frac{dS_p}{dt}&=\lambda_{s_p} R_s-k_2\left(1-e^{-\lambda S_p}\right)D S_p-\eta_{sp}S_p,\hspace{1.42cm}
			\frac{dS}{dt}	=\beta_1 P_g- \beta_2 S,\\
			\frac{dD}{dt}	&=\beta_2 S -k_1e^{-\lambda S_p} D-k_2\left(1-e^{-\lambda S_p}\right)D S_p,\hspace{0.95cm}
			\frac{dV}{dt}	= k_2\left(1-e^{-\lambda S_p}\right)D S_p-\delta_v V.
		\end{split}
		\label{simplified model}
		\right\}
	\end{equation}
We also perform the GSA of the simplified model \eqref{simplified model}. Due to the simplification, the PRCC values of parameters obtained for this model change from those of the full model \eqref{full model}, but the changes are not significant. However, the most positively and the most negatively sensitive parameters for each compartment remain  same in both models. Consequently, the results discussed in Sub-section \ref{Result: Scatter plots} remain unaffected in spite of simplification. Therefore, depending on the requirements, the simplified model \eqref{simplified model} may be used in place of the full model \eqref{full model}.

	\end{document}